\documentclass{amsart}
\usepackage[utf8]{inputenc}
\usepackage{amsmath}
\usepackage{amssymb}
\usepackage{wrapfig}
\usepackage{array}
\usepackage{subcaption}
\usepackage{graphicx}
\usepackage{biblatex} 
\usepackage[justification=centering]{caption}
\newtheorem{theorem}{Theorem}[section]
\newtheorem{proposition}[theorem]{Proposition}
\newtheorem{definition}{Definition}
\newtheorem{lemma}[theorem]{Lemma}
\newtheorem{corollary}[theorem]{Corollary}
\newtheorem*{hypothesis}{Hypothesis (F)}
\theoremstyle{definition}
\newtheorem{example}[theorem]{Example}
\theoremstyle{remark}
\newtheorem{remark}[theorem]{Remark}
\theoremstyle{plain}

\title[Discrete Distributions in the Surreal Numbers]{Discrete Distributions and Their Limit Theorems in the Surreal Numbers}
\author{John Harris}
\address{The University of Southern Mississippi\\
Long Beach, Mississippi}
\email{john.m.harris@usm.edu}
\date{November 6, 2026}

\begin{document}

\begin{abstract}
We extend the discrete uniform, hypergeometric and binomial
distributions to the surreal numbers, allowing infinitely many outcomes,
chips and trials. A uniform distribution on $\mathbb{N}$ gives each
outcome probability $1/\omega$; on the eventually periodic sets we
define the needed sizes explicitly, obtaining a probability that is
regular and strictly monotone. For an urn with infinitely many chips
and finitely many draws, the hypergeometric distribution is within an
infinitesimal of the binomial. With infinitely many trials, and under a
hypothesis on factorials of omnific integers that we show can be
satisfied, we prove a surreal de~Moivre--Laplace theorem, and we show
that the hypergeometric distribution with infinitely many draws is
locally normal with a finite-population correction. Near its center, it
is relatively close to the binomial distribution if and only if the
sampling fraction is infinitesimal, an exact form of a familiar rule of
thumb. For rare
successes we prove, without any hypothesis on factorials, a surreal
Poisson limit theorem and describe urns with finitely many red chips,
where the binomial approximation fails beyond the first success.
\end{abstract}

\maketitle


\section{Introduction}

Suppose that a man falls asleep under a tree, and the duration of his nap in minutes is a discrete random variable $X$. Any minute is as likely to be the minute he awakens as any other, and there is no upper bound on the duration. Suppose then that the space of $X$ is $\mathbb{N}$.

Of course, this doesn't make sense with a real-valued probability mass function, since there is no real number $c$ such that $\sum_{i=1}^\infty c = 1$. But with a surreal-valued pmf, it can make sense: we can let $c = \varepsilon = 1/\omega$. Then the pmf is $f(k) = \varepsilon$ for all $k \in \mathbb{N}$.

To assign probabilities to infinite events, such as the event that $X$
is even, we need a size for infinite subsets of $\mathbb{N}$. In
Section~2 we use the magnum of Lynch and Mackey \cite{LynchMackey},
which we define explicitly on the eventually periodic sets; the
resulting probability is regular and strictly monotone. In Section~3
we draw finitely many chips from an urn containing infinitely many,
and show that the hypergeometric distribution is within an
infinitesimal of the binomial. Section~4 lets the number of trials be
infinite. Under a hypothesis on factorials of omnific integers, which
we show can be satisfied, we prove a surreal de~Moivre--Laplace theorem
(Theorem~\ref{thm:dML}). Section~5 draws infinitely many chips: the
hypergeometric distribution is locally normal with a finite-population
correction (Theorem~\ref{thm:hyper-local}), and near its center it is
relatively close to the binomial distribution if and only if the
sampling fraction is infinitesimal (Theorem~\ref{thm:dichotomy}).
Section~6 treats rare successes, with a surreal Poisson limit theorem
(Theorem~\ref{thm:poisson}) and a description of urns with finitely
many red chips; no hypothesis on factorials is needed there. Section~7
collects open problems.

This work belongs to the study of non-Archimedean probability, which
Benci, Horsten and Wenmackers developed with hyperreal values
\cite{BHW2013,BHW2018}. The surreal numbers provide a setting in which
$\omega$, $\omega/2$, $\sqrt\omega$ and their reciprocals are all
available at once, together with an exponential function that lets us
transfer inequalities from the real numbers; and on the events we
consider, our probabilities are given by explicit formulas.

Infinitely many trials are also familiar from nonstandard analysis. If
$n$ is a hyperfinite integer, the binomial distribution with $n$ trials
is an internal finite probability space, its total mass is $1$ by
transfer, and theorems of de~Moivre--Laplace type hold there; see
Nelson \cite{Nelson}, especially Chapter~18. Loeb's construction
\cite{Loeb} then converts such hyperfinite spaces into standard,
real-valued measure spaces. Our setting differs in one essential
respect: there is no transfer principle from $\mathbb{Z}$ to the omnific
integers. The theorem of van den Dries and Ehrlich used in Section~4
transfers statements about the real exponential field, but the omnific
integers are not elementarily equivalent to $\mathbb{Z}$; for instance,
$(\sqrt2\,\omega)^{2}=2\omega^{2}$, so the equation $x^{2}=2y^{2}$ has
solutions in nonzero omnific integers. Consequently, factorials of
infinite omnific integers are not given to us, and we isolate what we
need as Hypothesis~(F) in Section~4; likewise, the normalization
$\sum_{k}b(k)=1$ is not automatic (Remark~\ref{rem:normalization}). In
return, $\omega$ is canonical rather than an arbitrary hyperfinite
integer, no ultrafilter has to be chosen, and many of our
probabilities, such as $2^{-\omega}$ and $\omega^{-1/2}$, are specific
surreal numbers.

\section{Extending the Uniform Distribution}

What if we want the probability that $X$ is even? Cardinality isn't up
to the task, since the even numbers and $\mathbb{N}$ have the same
cardinality. Natural density \cite[Chapter~III.1]{Tenenbaum}, the limit
of $|A\cap\{1,\dots,n\}|/n$ as $n\to\infty$, gives the expected answer
$\frac12$, but
it is too coarse in another way: it gives every finite set, and in
particular every single outcome, probability $0$, and it cannot tell
the odd numbers $\{1,3,5,\dots\}$ from $\{3,5,7,\dots\}$. We want a
surreal-valued size of sets that extends cardinality, is additive, and
gives a proper subset a strictly smaller size. On a natural class of
subsets of $\mathbb{N}$ such a size can be written down explicitly.

Write $\kappa_{A}(n)=|A\cap\{1,\dots,n\}|$ for the counting function of
$A\subseteq\mathbb{N}$. A set $P\subseteq\mathbb{N}$ is \emph{periodic}
if it is a union of residue classes modulo some $q\ge1$; if it consists
of $p$ such classes, put $d(P)=p/q$, its \emph{density}. (This does not depend
on the choice of $q$: passing from $q$ to a multiple $q'=qr$ replaces
each class by $r$ classes.)

\begin{definition}\label{def:eventually-periodic}
A set $A\subseteq\mathbb{N}$ is \emph{eventually periodic} if there is a
periodic set $P$ such that the symmetric difference $A\,\triangle\,P$ is
finite. We write $\mathcal{E}$ for the family of eventually periodic
sets.
\end{definition}

Note that the periodic set $P$ in Definition~\ref{def:eventually-periodic} is
unique: suppose $P$ and $P'$ both work. If $x\in P\setminus P'$, then
either $x\notin A$, so $x\in P\setminus A$, or $x\in A$, so
$x\in A\setminus P'$; the case $x\in P'\setminus P$ is symmetric.
Hence
\[
P\,\triangle\,P'\subseteq(P\,\triangle\,A)\cup(A\,\triangle\,P'),
\]
so $P\,\triangle\,P'$ is finite. It is also periodic: if $P$ and $P'$
are unions of residue classes modulo $q$ and $q'$, both are unions of
residue classes modulo $qq'$, and hence so is $P\,\triangle\,P'$.
Finally, a nonempty periodic set contains a whole residue class, which
is infinite; so the finite periodic set $P\,\triangle\,P'$ is empty,
and $P=P'$. We call $P$ the \emph{periodic part} of $A$.

\begin{definition}\label{def:magnum}
The \emph{magnum} of $A\in\mathcal{E}$, with periodic part $P$, is
\[
m(A)=d(P)\,\omega+|A\setminus P|-|P\setminus A|.
\]
\end{definition}

Since $A\setminus P$ and $P\setminus A$ are finite, $m(A)$ is a rational
multiple of $\omega$ plus an integer, and in particular an omnific
integer (see Section~4). For example, $m(\mathbb{N})=\omega$, the even
numbers and the odd numbers both have magnum $\omega/2$, the multiples
of $k$ have magnum $\omega/k$, $m(\{3,5,7,\dots\})=\omega/2-1$, and a
finite set has magnum equal to its number of elements.

The correction term can also be read off from counting functions. Since
\[
\kappa_{A}(n)-\kappa_{P}(n)
=\bigl|(A\setminus P)\cap\{1,\dots,n\}\bigr|-\bigl|(P\setminus A)\cap\{1,\dots,n\}\bigr|
\]
is eventually equal to $|A\setminus P|-|P\setminus A|$, we have
\begin{equation}\label{eq:magnum-limit}
m(A)=d(P)\,\omega+\lim_{n\to\infty}\bigl(\kappa_{A}(n)-\kappa_{P}(n)\bigr).
\end{equation}
Since $\kappa_{P}(n)$ is approximately $d(P)\,n$, this lets us read
$m(A)$ as the counting function of $A$ ``evaluated at $\omega$'', in
the spirit of Lynch and Mackey (Remark~\ref{rem:magnum-facts}). We
compare $\kappa_{A}$ with $\kappa_{P}$ rather than with $d(P)\,n$
because $\kappa_{P}(n)-d(P)\,n$ need not converge; for the even
numbers it alternates between $0$ and $-\frac12$.

\begin{proposition}\label{prop:magnum}
The family $\mathcal{E}$ is a Boolean algebra of subsets of
$\mathbb{N}$ containing every finite set and every residue class, and
the magnum $m\colon\mathcal{E}\to\mathbf{No}$ has the following
properties:
\begin{enumerate}
\item[(M1)] $m(A)=|A|$ for finite $A$;
\item[(M2)] $m(A\uplus B)=m(A)+m(B)$ for disjoint $A,B\in\mathcal{E}$;
\item[(M3)] $m(\mathbb{N})=\omega$;
\item[(M4)] for each $q\ge1$, the $q$ residue classes modulo $q$ have
equal magnums;
\item[(M5)] if $A\subsetneq B$, then $m(A)<m(B)$.
\end{enumerate}
Moreover, $m$ is the only function on $\mathcal{E}$ satisfying
(M1)--(M4).
\end{proposition}

\begin{proof}
If $A,B\in\mathcal{E}$ have periodic parts $P$ and $Q$, then
$(\mathbb{N}\setminus A)\triangle(\mathbb{N}\setminus P)=A\triangle P$
and $(A\cup B)\triangle(P\cup Q)\subseteq(A\triangle P)\cup(B\triangle Q)$
are finite, so $\mathbb{N}\setminus A$ and $A\cup B$ lie in $\mathcal{E}$,
with periodic parts $\mathbb{N}\setminus P$ and $P\cup Q$. Finite sets
(with $P=\varnothing$) and residue classes (with $P$ the class itself)
lie in $\mathcal{E}$.

(M1) For finite $A$, $P=\varnothing$, so $m(A)=|A|$. (M3) and
(M4) are immediate from the definition. (M2) If $A$ and $B$ are
disjoint, then $P\cap Q$ is periodic and contained in the finite set
$(A\triangle P)\cup(B\triangle Q)$, so $P\cap Q=\varnothing$. Hence
the periodic part of $A\uplus B$ is $P\uplus Q$, with
$d(P\uplus Q)=d(P)+d(Q)$. Since counting functions of disjoint sets
add, $\kappa_{A\uplus B}-\kappa_{P\uplus Q}
=(\kappa_{A}-\kappa_{P})+(\kappa_{B}-\kappa_{Q})$, so by
\eqref{eq:magnum-limit} the correction terms add as well. (M5) Write $B=A\uplus C$ with
$C=B\setminus A\in\mathcal{E}$ nonempty. If $C$ is finite, then
$m(C)=|C|\ge1$; if $C$ is infinite, its periodic part is nonempty, so
$d>0$ and $m(C)$ is a positive infinite number. Either way
$m(B)=m(A)+m(C)>m(A)$ by (M2).

For uniqueness, let $m'$ satisfy (M1)--(M4). By (M2)--(M4), each
residue class modulo $q$ has $m'$-value $\omega/q$, so
$m'(P)=d(P)\,\omega$ for periodic $P$. For $A\in\mathcal{E}$ with
periodic part $P$, the set $A\cup P$ is both $A\uplus(P\setminus A)$
and $P\uplus(A\setminus P)$, with $P\setminus A$ and $A\setminus P$
finite, so (M1) and (M2) give
$m'(A)=d(P)\,\omega+|A\setminus P|-|P\setminus A|=m(A)$.
\end{proof}

Since a finite set has magnum equal to its number of elements, the
following definition extends the usual definition of combinatorial
probability.

\begin{definition}\label{def:events}
Let $S\in\mathcal{E}$ be a nonempty outcome space. An \emph{event} is a
set $A\in\mathcal{E}$ with $A\subseteq S$. If each outcome in $S$ is
equally likely, the \emph{probability} that $A$ occurs is
$\dfrac{m(A)}{m(S)}$.
\end{definition}

For the nap of the Introduction, $S=\mathbb{N}$ and $m(S)=\omega$. Each
single minute has probability $\frac1\omega=\varepsilon$, in agreement
with $f(k)=\varepsilon$. The probability that $X$ is even is
$\frac{\omega/2}{\omega}=\frac12$, as expected, while the probability
that $X$ is an odd number greater than $1$ is
$\frac{\omega/2-1}{\omega}=\frac12-\varepsilon$. By (M5), every nonempty
event has positive probability, and a proper subevent has strictly
smaller probability; natural density has neither property.

The first of these properties is usually called \emph{regularity}.
Regular probability functions with non-Archimedean values were
developed by Benci, Horsten and Wenmackers \cite{BHW2013}, using
hyperreal rather than surreal numbers; 
among other things they treat a fair lottery on $\mathbb{N}$, in which every ticket has the same
infinitesimal probability, much like our nap. See \cite{BHW2018} for a
broader discussion of infinitesimal probabilities and the objections
to them. Their probability is defined on every subset of $\mathbb{N}$,
but its values are obtained through a non-constructive limit operation
(the ideal on which it rests comes from Zorn's lemma), whereas on
$\mathcal{E}$ our assignment is explicit and, by
Proposition~\ref{prop:magnum}, unique.

\begin{remark}[Magnums in general]\label{rem:magnum-facts}
The name ``magnum'' is taken from Lynch and Mackey \cite{LynchMackey},
who define magnums for a much larger class of countable sets, by
extending the counting function $\kappa_{A}$ to surnatural arguments and
evaluating at $\omega$. On $\mathcal{E}$ their values agree with ours;
for instance, they too give the even and the odd numbers magnum
$\omega/2$. Their theory also reaches sets outside $\mathcal{E}$: the
perfect squares, for example, receive magnum $\sqrt{\omega}$, which
would make the probability that $X$ is a perfect square equal to
$\sqrt{\omega}/\omega=\omega^{-1/2}$. Such values rely on the more
general machinery of their paper, in particular an Axiom of Extension
governing how counting functions are extended to surnatural arguments;
as the authors note, such extensions need not be unique. Since none of
this is needed here, we work only on $\mathcal{E}$, where everything is
explicit and, by Proposition~\ref{prop:magnum}, forced by
(M1)--(M4). Every event considered in this paper lies in $\mathcal{E}$:
the parity events above, the sets of red and white chips in
Examples~\ref{ex:multiples-of-5} and~\ref{ex:five-red}, and the tail
sets $\{k,k+1,\dots\}$ in Remark~\ref{rem:geometric}. Extending
Definition~\ref{def:events} beyond $\mathcal{E}$ is a natural problem
for the general theory of magnums.
\end{remark}

\begin{definition}
For any nonnegative surreal number $\eta$, we say that $x$ and $y$ are \emph{$\eta$-close} if there is a finite number $c$ such that $|x-y|  \leq c\eta$.
\end{definition}

Note that $\eta$-close is an equivalence relation, since it exhibits  

\emph{Reflexivity:} $| x - x | = 0 \leq\eta$,

\emph{Symmetry:} $| y - x | = | x - y | \leq c\eta$, and 
 
\emph{Transitivity:} $| x - z | \leq | x - y | + | y - z | \leq c\eta + d\eta = (c + d)\eta$.

\begin{proposition}\label{prop:closeness} For any finite numbers $a$ and $b$ and infinitesimals $\eta_a$ and $\eta_b$, if $a'$ is $\eta_a$-close to $a$ and if $b'$ is $\eta_b$-close to $b$, then 
\begin{enumerate}
\item $a' \pm b'$ is  $\eta$-close to $a \pm b$,

\item $a'b'$ is $\eta$-close to $ab$, and  

\item  if $|b| > r$, for some $r \in \mathbb{R}^+$, $\frac{a'}{b'}$ is $\eta$-close to $\frac{a}{b}$,
\end{enumerate}
where $\eta = \max\{\eta_a,\eta_b\}$.
\end{proposition}

\begin{proof}
\begin{enumerate}
\item 
\begin{align*}
|(a' \pm b') - (a \pm b) | = |(a' - a) \pm (b' - b) | &\leq  |a' - a| + |b' - b | \\ &\leq c_a\eta_a + c_b\eta_b \\ &\leq (c_a + c_b)\eta.
\end{align*}

\item
Note that $a' = a \;+$ an infinitesimal, so $a'$ is also finite, and thus
\begin{align*}
|a'b' - ab| = |a'b' - a'b + a'b -ab| &= |a'(b' - b) + (a' - a)b| \\ &\leq |a'(b' - b)| + |(a' - a)b| \\
&\leq  |a'|c_b\eta_b + |b|c_a\eta_a \\ &\leq (|a'|c_b + |b|c_a)\eta.
\end{align*}

\item
 As above, $a'$ and $b'$ are both finite, and further, as $b' =  b \;+$ an infinitesimal, $|b'| \geq \frac{r}{2}$.  Hence, 

\begin{align*}
 \left|\frac{a'}{b'} - \frac{a}{b}\right| =  \left|\frac{a'b - ab'}{bb'}\right| =  \left|\frac{a'b - a'b' + a'b' - ab'}{bb'}\right| &= \frac{|a'(b - b') + (a' - a)b'|}{|bb'|} \\
&\leq \frac{|a'(b - b')|+ |(a' - a)b'|}{|bb'|} \\
&\leq \frac{|a'|c_b\eta_b + |b'|c_a\eta_a}{|bb'|} \\ & \leq \frac{|a'|c_b + |b'|c_a}{|bb'|}\eta \\ & \leq \frac{2(|a'|c_b + |b'|c_a)}{r^2}\eta
\end{align*}
\end{enumerate}
\end{proof}
\section{Extending the Hypergeometric Distribution}

Suppose that we have an urn containing $a\omega$ red chips and $b\omega$ white chips, 
for some positive reals $a$ and $b$, and we draw $n$ chips from the urn without replacement.
What is the probability that exactly $k$ chips are red? 
(For real $a>0$, the normal form of $a\omega$ is $a\omega^{1}$, so $a\omega$ is
an omnific integer; see Section~4.)

Here, for finite $k$ and any surreal number $x$, we put
\begin{equation}\label{eq:binom-finite}
\binom xk=\frac{x(x-1)\cdots(x-k+1)}{k!},
\end{equation}
a polynomial in $x$, so the binomial coefficients below make sense.

\begin{align*}
\frac{\binom{a\omega}{k} \binom{b\omega}{n-k}}{\binom{(a + b) \omega}{n}} &= \frac{(a\omega)(a\omega - 1) \cdots (a\omega - k + 1)(b\omega)(b\omega - 1) \cdots (b\omega - n + k + 1) n!}{k! (n-k)! ((a + b) \omega) ((a + b) \omega - 1)\cdots ((a+b)\omega - n  + 1)} \\
 &=\binom{n}{k} \frac{(a\omega)(a\omega - 1) \cdots (a\omega - k + 1)(b\omega)(b\omega - 1) \cdots (b\omega - n + k + 1)}{((a + b) \omega) ((a + b) \omega - 1)\cdots ((a+b)\omega - n  + 1)} \\
 &=\binom{n}{k} \frac{a\omega}{(a + b) \omega} \cdots \frac{a\omega - k + 1}{(a + b) \omega - k + 1} \frac{b\omega}{(a+b)\omega - k} \cdots \frac{b\omega - n + k + 1}{(a+b)\omega - n  + 1}
\end{align*}
 \medskip

Note that the first fraction on the right is $\frac{a}{a + b}$, 
while the second can be written as $\frac{a - \varepsilon}{a + b - \varepsilon}$, 
which is $\varepsilon$-close to $\frac{a}{a + b}$.
In fact, the first $k$ fractions on the right are all $\varepsilon$-close to $\frac{a}{a + b}$, 
while the remaining $n-k$ fractions are $\varepsilon$-close to $\frac{b}{a + b}$.
Thus, the product of all of the fractions is 
$\varepsilon$-close to $\left(\frac{a}{a + b}\right)^k \left(\frac{b}{a + b}\right)^{n-k}$, 
and so, the probability that exactly $k$ red chips are drawn is $\varepsilon$-close to 
$\binom{n}{k} \left(\frac{a}{a + b}\right)^k \left(\frac{b}{a + b}\right)^{n-k}$.

Also, note that if we replace $a\omega$ and $b\omega$ 
with $a\omega + c$ and $b\omega + d$ in our example, for any $c, d \in \mathbb{Z}$, 
the probability is still $\varepsilon$-close to $\binom{n}{k} \left(\frac{a}{a + b}\right)^k \left(\frac{b}{a + b}\right)^{n-k}$.

Allowing $a$ and $b$ to be arbitrary positive reals costs nothing in the
argument above. If we want the urn to be a concrete set of chips,
numbered by the natural numbers, with the red and the white chips
forming eventually periodic sets, then $a$ and $b$ are rational, since
the magnum of an eventually periodic set is a rational multiple of
$\omega$ plus an integer. (Outside $\mathcal{E}$, in the broader theory
of Lynch and Mackey, other values occur; for example, the perfect
squares have magnum $\sqrt\omega$. See Remark~\ref{rem:magnum-facts}.)

\begin{example}\label{ex:multiples-of-5}
Number the chips by the natural numbers, and let the chips numbered
$5,10,15,\dots$ be red and the rest white. The red chips have magnum
$\omega/5$, just as the even numbers have magnum $\omega/2$, and the
white chips have magnum $\omega-\omega/5=4\omega/5$. So $a=\frac15$ and
$b=\frac45$, and if we draw $n$ chips, the probability that exactly $k$
are red is $\varepsilon$-close to
\[
\binom nk\Bigl(\frac15\Bigr)^{k}\Bigl(\frac45\Bigr)^{n-k}.
\]
For instance, the probability that a single draw is red is exactly
$\frac{\omega/5}{\omega}=\frac15$, matching the intuition that one
natural number in five is a multiple of $5$.
\end{example}

If instead only finitely many chips are red, say those numbered
$1,2,3,4,5$, then the red chips have magnum $5$ and the white chips
$\omega-5$. This is formally the case $a=0$, $b=1$, $c=5$, $d=-5$ of
the remark above. The same computation works when $a=0$, but it only
tells us that drawing any red chip at all has probability
$\varepsilon$-close to $0$. Section~\ref{sec:rare}
gives a much finer description of this urn
(Example~\ref{ex:five-red}).

\begin{definition}
\label{def:uniform-close}For any nonnegative surreal number $\eta$, we say that $X$ and $Y$ are \emph{uniformly $\eta$-close} 
if there is a finite number $c$ such that $|f_X(x) - f_Y(x)|\leq c\eta$.
\end{definition}
\begin{theorem}\label{thm:uniform-close-exp}
If $X$ is \textit{uniformly $\eta$-close} to $Y$ and if $g(x)$ is absolutely summable on the space of $X$, 
then $E[g(X)]$ is \textit{$\eta$-close} to $E[g(Y)]$.
\end{theorem}
\begin{proof}
\begin{align*}| E[g(X)] - E[g(Y)] | &= \left|\sum_{x} g(x)f_X(x) - \sum_{x} g(x)f_Y(x)\right| \\ &= \left|\sum_{x} g(x)f_X(x) - g(x)f_Y(x)\right| \\ & = \left|\sum_{x} g(x)( f_X(x) - f_Y(x) )\right| \\ & \leq \sum_{x} |g(x)( f_X(x) - f_Y(x) )| \\ & \leq \sum_{x} |g(x)| | f_X(x) - f_Y(x) | \\ &\leq \sum_{x} |g(x)| c\eta \\ & =\left( \sum_{x} |g(x)| \right) c\eta \\ & \leq d\eta.
\end{align*}
\end{proof}
Thus, in our urn example, since our distribution is uniformly $\varepsilon$-close to a binomial distribution, 
our mean is $\varepsilon$-close to the binomial mean $n\frac{a}{a+b}$.  

\section{Infinitely Many Trials}

In Section 3, the number of draws was finite. Here, we let the number of
trials itself be infinite. Two new issues appear at once. First, the
number of successes in $\omega$ trials need not be an ordinary integer:
it ranges over the \emph{omnific integers} $\mathbf{Oz}$, the surreals
$x$ with $x=\{x-1\mid x+1\}$ \cite{Conway}. Between $0$ and $\omega$ lie
not only $0,1,2,\dots,\omega-1,\omega$ but also $\omega/2$,
$\omega^{1/2}$, $\pi\omega$, $\omega-\omega^{1/2}$, and uncountably many
others. Second, individual probabilities near the center of the
distribution are infinitesimal, so the absolute notion of
$\eta$-closeness from Section 2 is too coarse, and we need a relative
version.

\subsection{Tools}

We use Gonshor's exponential $\exp$ on $\mathbf{No}$ and its inverse
$\log$ on the positive surreals \cite{Gonshor}, and for $t>0$ we set
$t^{s}=\exp(s\log t)$. The key fact is the theorem of van den Dries and
Ehrlich \cite{vdDE} that $(\mathbf{No},+,\cdot,<,\exp)$ is an elementary
extension of $(\mathbb{R},+,\cdot,<,\exp)$. Consequently any inequality
of the form ``for all $t$ with $|t|\le r$, $|F(t)|\le C|t|^m$'', where
$F$ is built from field operations, $\exp$ and $\log$ and $r,C$ are
real, holds in $\mathbf{No}$ as soon as it holds in $\mathbb{R}$. We
call this \emph{transfer}. For instance, for all surreal $t$,
\begin{equation}\label{eq:transfer}
|t|\le 1 \implies |e^{t}-1|\le 2|t|,
\qquad
|t|\le \tfrac12 \implies \bigl|\log(1+t)-t\bigr|\le t^{2}.
\end{equation}

\begin{definition}
Let $\eta\ge 0$. Nonzero surreals $x$ and $y$ are \emph{relatively
$\eta$-close} if $x/y$ is $\eta$-close to $1$.
\end{definition}

Note that if $y$ is finite and $x$ and $y$ are relatively $\eta$-close, 
then $x$ and $y$ are $\eta$-close. Similarly, if $y$ is not infinitesimal 
and $x$ and $y$ are $\eta$-close, $x$ and $y$ are relatively $\eta$-close.

\begin{lemma}\label{lem:relative}
Let $\eta$ be infinitesimal. Relative $\eta$-closeness is an
equivalence relation on the nonzero surreals. Moreover, if $x$ is
relatively $\eta$-close to $x'$ and $y$ is relatively $\eta$-close to
$y'$, then $xy$ is relatively $\eta$-close to $x'y'$, and $x/y$ is
relatively $\eta$-close to $x'/y'$.
\end{lemma}

\begin{proof}
The numbers $x$ and $y$ may be infinite or infinitesimal, so we do not
apply Proposition~\ref{prop:closeness} to them directly. Instead we
apply it to their ratios, which are $\eta$-close to $1$. In every case
below the reference values are $a=b=1$, which are finite with
$|b|>\tfrac12$, so the hypotheses of Proposition~\ref{prop:closeness}
hold. Note also that $1$ is trivially $\eta$-close to $1$.

\emph{Reflexivity:} $x/x=1$.

\emph{Symmetry:} Suppose $r=x/y$ is $\eta$-close to $1$. By part~(3)
with $a'=a=1$, $b'=r$ and $b=1$, the number $y/x=1/r$ is $\eta$-close
to $1/1=1$.

\emph{Transitivity:} Suppose $r=x/y$ and $s=y/z$ are $\eta$-close to
$1$. By part~(2) with $a'=r$, $b'=s$ and $a=b=1$, the number $x/z=rs$
is $\eta$-close to $1$.

\emph{Products:} Let $r=x/x'$ and $s=y/y'$, both $\eta$-close to $1$.
Then $(xy)/(x'y')=rs$, which is $\eta$-close to $1$ by part~(2).

\emph{Quotients:} With $r$ and $s$ as above,
$(x/y)/(x'/y')=r/s$, which is $\eta$-close to $1$ by part~(3).
\end{proof}

By \eqref{eq:transfer}, if $u$ is
$\eta$-close to $0$ and $\eta$ is infinitesimal, then $e^{u}$ is
relatively $\eta$-close to $1$.

Throughout Sections~4 and~5, $0<p<1$ is real and $q=1-p$; we also put
\[
\varphi(x)=\frac{1}{\sqrt{2\pi}}e^{-x^{2}/2}.
\]

\subsection{Factorials of omnific integers}

For finite $k$, the binomial coefficient $\binom nk$ was defined in
Section 3 by \eqref{eq:binom-finite}, as a polynomial in $n$. This
already gives some interesting values.

\begin{example}
Flip a fair coin $\omega$ times. The probability of no heads is
$2^{-\omega}=\exp(-\omega\log 2)$, of exactly one head is
$\omega\,2^{-\omega}$, and, for finite $k$, of exactly $k$ heads is
$\binom{\omega}{k}2^{-\omega}$. By transfer, $2^{-\omega}<\omega^{-m}$
for every finite $m$, so each of these is smaller than every power
$\varepsilon^{m}$: a finite number of heads in $\omega$ fair flips is
``doubly'' negligible.
\end{example}

For infinite $k$ we need more. No canonical factorial on $\mathbf{Oz}$ is
known to us, so we isolate what we use as a hypothesis.

\begin{hypothesis}
There is a function $k\mapsto k!$ on the nonnegative omnific integers
such that
\begin{enumerate}
\item[(F1)] $0!=1$ and $(k+1)!=(k+1)\,k!$ for all $k$;
\item[(F2)] for every infinite $k$, $k!$ is relatively $(1/k)$-close to
$\sqrt{2\pi k}\,(k/e)^{k}$.
\end{enumerate}
Here $(k/e)^{k}$ means $\exp(k\log k-k)$.
\end{hypothesis}

\begin{remark}\label{rem:F-exists}
Functions satisfying (F) exist, though not canonically. (F1) only
links $k$ to $k\pm 1$, so it determines $k!$ separately on each coset
$k_{0}+\mathbb{Z}$. On the coset of $0$, use the ordinary factorial. On
each infinite coset, choose a representative $k_{0}$ (using global
choice, or working below a fixed bound), set
$k_{0}!=\sqrt{2\pi k_{0}}\,(k_{0}/e)^{k_{0}}$, and extend by (F1) in
both directions.

We check (F2). Write $S(k)=\sqrt{2\pi k}\,(k/e)^{k}$, so that
$\log S(k)=\bigl(k+\tfrac12\bigr)\log k-k+\tfrac12\log 2\pi$, and let
$k=k_{0}+j$ with $j\ge1$ finite. (The case $j\le-1$ is the same, using
$\log k_{0}!-\log k!=\sum_{i=0}^{|j|-1}\log(k_{0}-i)$.) Every number $t$
that we feed into $\log(1+t)$ below is infinitesimal, so the second
inequality in \eqref{eq:transfer} gives $|\log(1+t)-t|\le t^{2}$, and thus, $|\log(1+t)|\le2|t|$.
\begin{itemize}
\item \emph{The factorials.} By (F1),
\[
\log k!-\log k_{0}!=\sum_{i=1}^{j}\log(k_{0}+i)=j\log k_{0}+u_{1},
\qquad u_{1}=\sum_{i=1}^{j}\log\Bigl(1+\frac{i}{k_{0}}\Bigr),
\]
and $|u_{1}|\le\sum_{i=1}^{j}\frac{2i}{k_{0}}=\frac{j(j+1)}{k_{0}}$.
\item \emph{The Stirling expressions.} Put $t=j/k_{0}$, so that
$\log k=\log k_{0}+\log(1+t)$. The terms $-k$ and $-k_{0}$ contribute
$-j$, so
\[
\log S(k)-\log S(k_{0})
=\Bigl(k_{0}+j+\tfrac12\Bigr)\log k-\Bigl(k_{0}+\tfrac12\Bigr)\log k_{0}-j
=j\log k_{0}+u_{2},
\]
where $u_{2}=\bigl(k_{0}+j+\tfrac12\bigr)\log(1+t)-j$. Write
$\log(1+t)=t+R$ with $|R|\le t^{2}$. Since
$\bigl(k_{0}+j+\tfrac12\bigr)t=j+\frac{(j+\frac12)j}{k_{0}}$, and
$k_{0}+j+\tfrac12\le2k_{0}$,
\[
u_{2}=\frac{(j+\frac12)j}{k_{0}}+\Bigl(k_{0}+j+\tfrac12\Bigr)R,
\qquad
|u_{2}|\le\frac{(j+\frac12)j}{k_{0}}+2k_{0}\cdot\frac{j^{2}}{k_{0}^{2}}
=\frac{(j+\frac12)j+2j^{2}}{k_{0}}.
\]
\item \emph{Comparison.} Since $k_{0}!=S(k_{0})$, subtracting the two
displays gives $\log k!-\log S(k)=u_{1}-u_{2}$, and
$|u_{1}-u_{2}|\le c_{j}/k_{0}$ with
$c_{j}=j(j+1)+(j+\frac12)j+2j^{2}$ finite. By the consequence of
\eqref{eq:transfer} noted after Lemma~\ref{lem:relative},
$k!/S(k)=e^{u_{1}-u_{2}}$ is relatively $(1/k_{0})$-close to $1$. Finally
$k\le2k_{0}$, so $1/k_{0}\le2/k$, and (F2) holds at $k$, with a constant
depending on $j$.
\end{itemize}
(A third-order expansion of $\log(1+t)$ shows that the two differences
in fact agree more closely: both equal
$j\log k_{0}+\frac{j(j+1)}{2k_{0}}+O(j^{3}/k_{0}^{2})$, so $k!/S(k)$ is
even relatively $(1/k_{0}^{2})$-close to $1$. We do not need this.)

Different choices change
$k!$ by a factor relatively $(1/k)$-close to $1$, so every statement
below holds for every choice. Whether there is a
natural choice, for example via a surreal Gamma function, is open.

One would expect $k!$ to be an omnific integer. Recall that a surreal
number is an omnific integer exactly when every exponent in its normal
form is nonnegative and the coefficient of $\omega^{0}$ is an integer
\cite{Conway}. For the cosets that occur in this paper, the
construction above already produces omnific integers.

\emph{Claim.} Let $k_{0}=c\,\omega^{x}$ with $c$ and $x$ positive reals.
Then $S(k_{0})=C\,\omega^{Y}$ for a positive real $C$ and an infinite
$Y$, and if $k_{0}!=S(k_{0})$, then $k!$ is an omnific integer for
every $k$ in $k_{0}+\mathbb{Z}$.

\emph{Proof of the claim.} Call a surreal number \emph{purely infinite}
if every exponent in its normal form is positive. By Gonshor's theory
\cite{Gonshor}, $\exp$ maps the purely infinite numbers onto the
monomials $\omega^{y}$, and so $\log$ maps monomials to purely infinite
numbers. Now
\[
k_{0}\log k_{0}-k_{0}
=(c\log c-c)\,\omega^{x}+c\,\omega^{x}\log\bigl(\omega^{x}\bigr).
\]
The first term is purely infinite since $x>0$. In the second,
$\log(\omega^{x})$ is purely infinite, and multiplying by $\omega^{x}$
raises every exponent by $x$, so it too is purely infinite. Hence
$\exp(k_{0}\log k_{0}-k_{0})=\omega^{y}$ for some $y$, and since
$\sqrt{2\pi k_{0}}=\sqrt{2\pi c}\,\omega^{x/2}$,
\[
S(k_{0})=\sqrt{2\pi c}\;\omega^{Y},\qquad Y=\tfrac{x}{2}+y.
\]
To see that $Y$ is infinite, fix a finite $m$. The real inequality
$\sqrt{2\pi t}\,(t/e)^{t}\ge t^{m}$ holds for all sufficiently large
real $t$, so by transfer $S(k_{0})\ge k_{0}^{m}=c^{m}\omega^{mx}$. If
$Y<mx$, then $\omega^{Y}/\omega^{mx}$ would be infinitesimal, which is
impossible; so $Y\ge mx$ for every finite $m$, and $Y$ is infinite.

In particular $S(k_{0})$ is an omnific integer, and going upward,
$(k_{0}+j)!=S(k_{0})\,(k_{0}+1)\cdots(k_{0}+j)$ is an omnific integer
for every finite $j\ge0$, because the omnific integers form a ring.
Going downward, for finite $j\ge1$,
\[
(k_{0}-j)!=\frac{S(k_{0})}{k_{0}(k_{0}-1)\cdots(k_{0}-j+1)}
=\frac{\sqrt{2\pi c}\;\omega^{Y}}{c^{j}\,\omega^{jx}}
\prod_{i=0}^{j-1}\Bigl(1-\frac{i}{k_{0}}\Bigr)^{-1}.
\]
Each factor $(1-i/k_{0})^{-1}=\sum_{n\ge0}(i/c)^{n}\omega^{-nx}$ is a
power series in $\omega^{-x}$ with real coefficients, and so is their
product. Therefore $(k_{0}-j)!$ is a sum of real multiples of the
monomials $\omega^{Y-(j+n)x}$, $n\ge0$. Since $Y$ is infinite and
$(j+n)x$ is a real number, every exponent $Y-(j+n)x$ is positive. So
there are no terms with exponent $\le0$, and $(k_{0}-j)!$ is an omnific
integer. This proves the claim.

The claim covers $k_{0}=\omega$, and more generally $k_{0}=a\omega$ and
$k_{0}=a\omega^{1/2}$ for every positive real $a$; these include the
centers $np$ of the examples in Section~5. It does not cover every coset
that our theorems range over: Theorem~\ref{thm:dML}, for instance,
involves $k=np+x\sqrt{npq}$, such as $k=\omega/2+\omega^{1/2}$ when
$n=\omega$, $p=\frac12$ and $x=2$. For a general infinite omnific integer
$k_{0}$, such as $\omega+\omega^{1/2}$, the number $S(k_{0})$ has more
than one term,
and we do not know whether it, or its quotients by
$k_{0}(k_{0}-1)\cdots(k_{0}-j+1)$, are omnific integers. Going upward,
this can be repaired in any case. Write $S=S(k_{0})$. There is an omnific integer $\lfloor S\rfloor$ with
$\lfloor S\rfloor\le S<\lfloor S\rfloor+1$ \cite{Conway}, and thus,
$S-1<\lfloor S\rfloor\le S$. Dividing by $S>0$,
\[
1-\frac1S<\frac{\lfloor S\rfloor}{S}\le1,
\qquad\text{so}\qquad
\Bigl|\frac{\lfloor S\rfloor}{S}-1\Bigr|<1\cdot\frac1S,
\]
which says that $\lfloor S\rfloor/S$ is $(1/S)$-close to $1$, with
constant $1$. Moreover $S\ge k_{0}$: the real inequality
$\sqrt{2\pi x}\,(x/e)^{x}\ge x$ holds for all $x\ge3$ and is of the kind
covered by transfer, and $k_{0}>3$. Hence $1/S\le1/k_{0}$, and
$\lfloor S\rfloor/S$ is also $(1/k_{0})$-close to $1$. So we may set
$k_{0}!=\lfloor S\rfloor$ instead: this multiplies every $k!$ in the
coset by $\lfloor S\rfloor/S$, which by Lemma~\ref{lem:relative} does
not affect (F2). Then
$(k_{0}+j)!=k_{0}!\,(k_{0}+1)\cdots(k_{0}+j)$ is an omnific integer for
every finite $j\ge0$, because the omnific integers form a ring. Going
downward, however,
\[
(k_{0}-j)!=\frac{k_{0}!}{k_{0}(k_{0}-1)\cdots(k_{0}-j+1)},
\]
so for all of these to be omnific integers, $k_{0}!$ must be divisible
by $k_{0}(k_{0}-1)\cdots(k_{0}-j+1)$ for every finite $j$ simultaneously.
For general $k_{0}$ we do not know whether such an omnific integer,
relatively $(1/k_{0})$-close to $S(k_{0})$, exists. Nothing below
depends on $k!$ being an omnific integer: we use only (F1) and (F2).
\end{remark}

Assuming (F), set $\binom nk=\dfrac{n!}{k!\,(n-k)!}$ for omnific
$0\le k\le n$. By (F1) this agrees with \eqref{eq:binom-finite} when $k$
is finite.

The size of $\binom nk$ is governed by the \emph{binary entropy function}
\[
\mathrm{H}(t)=-t\log t-(1-t)\log(1-t),\qquad 0<t<1,
\]
the Shannon entropy of a single trial with success probability $t$
\cite{CoverThomas}.

\begin{lemma}\label{lem:stirling-binom}
Assume (F). Let $m,j$ be omnific with $j$ and $m-j$ infinite, and put
$t=j/m$. Then $\binom mj$ is relatively $\frac{1}{\min(j,\,m-j)}$-close to
\[
\frac{e^{m\mathrm{H}(t)}}{\sqrt{2\pi m\,t(1-t)}}.
\]
\end{lemma}

\begin{proof}
Write $S(k)=\sqrt{2\pi k}\,(k/e)^{k}$ for the Stirling expression in (F2),
and put
\[
\eta=\frac{1}{\min(j,\,m-j)}=\max\Bigl(\frac1j,\,\frac1{m-j}\Bigr).
\]
By (F2), $m!$, $j!$ and $(m-j)!$ are
relatively close to $S(m)$, $S(j)$ and $S(m-j)$, with error terms
$\frac1m$, $\frac1j$ and $\frac1{m-j}$. Each of these is at most $\eta$:
for $\frac1j$ and $\frac1{m-j}$ this is the definition of $\eta$, and since $0<j<m$ we have
$\frac1m<\frac1j\le\eta$. So all three are relatively $\eta$-close, and
by Lemma~\ref{lem:relative},
\[
\binom mj=\frac{m!}{j!\,(m-j)!}
\ \text{is relatively $\eta$-close to}\
\frac{S(m)}{S(j)\,S(m-j)}.
\]
It remains to compute the right-hand side. Since
$\log S(k)=\bigl(k+\tfrac12\bigr)\log k-k+\tfrac12\log 2\pi$, and the
terms $-m+j+(m-j)$ cancel,
\[
\log\frac{S(m)}{S(j)\,S(m-j)}
=\Bigl(m+\tfrac12\Bigr)\log m-\Bigl(j+\tfrac12\Bigr)\log j
-\Bigl(m-j+\tfrac12\Bigr)\log(m-j)-\tfrac12\log 2\pi.
\]
Now substitute $j=tm$ and $m-j=(1-t)m$, so that $\log j=\log m+\log t$
and $\log(m-j)=\log m+\log(1-t)$. The coefficient of $\log m$ becomes
\[
\Bigl(m+\tfrac12\Bigr)-\Bigl(j+\tfrac12\Bigr)-\Bigl(m-j+\tfrac12\Bigr)=-\tfrac12,
\]
and the remaining terms are
\begin{align*}
&-\Bigl(j+\tfrac12\Bigr)\log t-\Bigl(m-j+\tfrac12\Bigr)\log(1-t)\\
&\qquad=-mt\log t-m(1-t)\log(1-t)-\tfrac12\log t-\tfrac12\log(1-t)\\
&\qquad=m\mathrm{H}(t)-\tfrac12\log\bigl(t(1-t)\bigr).
\end{align*}
Altogether,
\begin{align*}
\log\frac{S(m)}{S(j)\,S(m-j)}
&=m\mathrm{H}(t)-\tfrac12\log m-\tfrac12\log\bigl(t(1-t)\bigr)-\tfrac12\log 2\pi\\
&=m\mathrm{H}(t)-\tfrac12\log\bigl(2\pi m\,t(1-t)\bigr),
\end{align*}
and exponentiating gives the expression in the lemma.
\end{proof}

\subsection{The binomial distribution with infinitely many trials}

Let $n$ be an infinite omnific integer. For omnific $0\le k\le n$ define
\[
b(k)=b(k;n,p)=\binom nk p^{k}q^{n-k}.
\]
We do not claim that $\sum_{k}b(k)=1$. The index set is uncountable,
and no notion of summation over it is available, so the binomial
theorem cannot simply be invoked (see Remark~\ref{rem:normalization}).
What we can do is describe $b$ locally, and the description is the
de~Moivre--Laplace theorem \cite{Feller}. In its proof, the binomial
weight $p^{k}q^{n-k}$ combines with $e^{n\mathrm{H}(k/n)}$ into a single
exponential, whose exponent is $-n$ times
\[
D(t\,\|\,p)=t\log\frac tp+(1-t)\log\frac{1-t}{q},\qquad 0<t<1.
\]
This is the relative entropy (Kullback--Leibler divergence) of a
Bernoulli$(t)$ distribution from a Bernoulli$(p)$ distribution
\cite{CoverThomas}. It is nonnegative and vanishes only at $t=p$; we
use only its Taylor expansion about $t=p$, which we now make precise.

\begin{lemma}\label{lem:D-taylor}
Let $0<p<1$ be real, $q=1-p$, and
$C_{p}=\frac23\bigl(\frac1{p^{2}}+\frac1{q^{2}}\bigr)$. Then for every
surreal $s$ with $|s|\le\frac12\min(p,q)$,
\[
\Bigl|D(p+s\,\|\,p)-\frac{s^{2}}{2pq}\Bigr|\le C_{p}|s|^{3}.
\]
\end{lemma}

\begin{proof}
First let $s$ be real. Put
\[
f(s)=D(p+s\,\|\,p)=(p+s)\log\frac{p+s}{p}+(q-s)\log\frac{q-s}{q},
\]
which is smooth for $-p<s<q$, and in particular on the closed interval
$|s|\le\frac12\min(p,q)$. Differentiating,
\begin{align*}
f'(s)&=\log\frac{p+s}{p}+1-\log\frac{q-s}{q}-1
=\log\Bigl(1+\frac sp\Bigr)-\log\Bigl(1-\frac sq\Bigr),\\
f''(s)&=\frac{1}{p+s}+\frac{1}{q-s},\\
f'''(s)&=-\frac{1}{(p+s)^{2}}+\frac{1}{(q-s)^{2}}.
\end{align*}
Thus $f(0)=0$, $f'(0)=0$ and $f''(0)=\frac1p+\frac1q=\frac{1}{pq}$,
since $p+q=1$. On the interval $|s|\le\frac12\min(p,q)$ we have
$p+s\ge\frac p2$ and $q-s\ge\frac q2$, so
\[
|f'''(s)|\le\frac{4}{p^{2}}+\frac{4}{q^{2}}.
\]
By Taylor's theorem with the Lagrange form of the remainder, for each
such $s$ there is $\xi$ between $0$ and $s$ with
\[
f(s)=f(0)+f'(0)s+\frac{f''(0)}{2}s^{2}+\frac{f'''(\xi)}{6}s^{3}
=\frac{s^{2}}{2pq}+\frac{f'''(\xi)}{6}s^{3},
\]
and therefore
\[
\Bigl|f(s)-\frac{s^{2}}{2pq}\Bigr|\le\frac16\Bigl(\frac{4}{p^{2}}+\frac{4}{q^{2}}\Bigr)|s|^{3}
=C_{p}|s|^{3}.
\]
This proves the lemma for real $s$. The statement ``for all $s$ with
$|s|\le\frac12\min(p,q)$,
$\bigl|D(p+s\,\|\,p)-\frac{s^{2}}{2pq}\bigr|\le C_{p}|s|^{3}$'' is an
inequality of the kind described in Section~4.1: $D$ is built from
field operations and $\log$, and the parameters $p$, $q$ and $C_{p}$
are real. By transfer it therefore holds for all surreal $s$ as well.
Note that we transfer only the \emph{conclusion} of Taylor's theorem,
a single inequality, and not Taylor's theorem itself.
\end{proof}

\begin{theorem}[Surreal de Moivre--Laplace]\label{thm:dML}
Assume (F). Let $n$ be an infinite omnific integer and let $k$ be an
omnific integer with $x=(k-np)/\sqrt{npq}$ finite. Then $b(k;n,p)$ is relatively
$n^{-1/2}$-close to
\[
\frac{\varphi(x)}{\sqrt{npq}}.
\]
\end{theorem}

\begin{proof}
Put $t=k/n$. Since $k=np+x\sqrt{npq}$, we have
\[
t=p+s,\qquad\text{where } s=\frac{x\sqrt{npq}}{n}=x\sqrt{pq}\cdot\frac{1}{\sqrt n}.
\]
Here $x$ is finite, $\sqrt{pq}$ is real, and $1/\sqrt n$ is infinitesimal
because $n$ is infinite, so $s$ is infinitesimal. In particular
$|s|\le\frac12\min(p,q)$, and therefore
\[
k=n(p+s)\ge\tfrac12 np
\qquad\text{and}\qquad
n-k=n(q-s)\ge\tfrac12 nq.
\]
Since $n$ is infinite and $p,q$ are positive reals, both $k$ and $n-k$
are infinite, so Lemma~\ref{lem:stirling-binom} applies with $m=n$ and
$j=k$. Its error term satisfies
\[
\frac{1}{\min(k,\,n-k)}\le\frac1k+\frac1{n-k}\le\frac{2}{np}+\frac{2}{nq}=\frac{2}{pq}\cdot\frac1n,
\]
a real multiple of $1/n$. So Lemma~\ref{lem:stirling-binom} gives
\begin{equation}\label{eq:dML-binom}
\binom nk\ \text{is relatively } n^{-1}\text{-close to}\
\frac{e^{n\mathrm{H}(t)}}{\sqrt{2\pi n\,t(1-t)}}.
\end{equation}

\emph{Step 1: the exponent.} Since $k=nt$ and $n-k=n(1-t)$,
\[
p^{k}q^{n-k}=\exp\bigl(k\log p+(n-k)\log q\bigr)
=\exp\bigl(n\bigl[t\log p+(1-t)\log q\bigr]\bigr).
\]
Multiplying both sides of \eqref{eq:dML-binom} by this exact factor does
not change their ratio, so $b(k)=\binom nk p^{k}q^{n-k}$ is relatively
$n^{-1}$-close to
\[
\frac{\exp\bigl(n\bigl[\mathrm{H}(t)+t\log p+(1-t)\log q\bigr]\bigr)}{\sqrt{2\pi n\,t(1-t)}}.
\]
The bracket in the exponent is $-D(t\,\|\,p)$:
\begin{align*}
\mathrm{H}(t)+t\log p+(1-t)\log q
&=-t\log t-(1-t)\log(1-t)+t\log p+(1-t)\log q\\
&=-\Bigl[t\log\frac tp+(1-t)\log\frac{1-t}{q}\Bigr]
=-D(t\,\|\,p).
\end{align*}
Hence
\begin{equation}\label{eq:dML-D}
b(k)\ \text{is relatively } n^{-1}\text{-close to}\
\frac{e^{-nD(t\|p)}}{\sqrt{2\pi n\,t(1-t)}}.
\end{equation}

\emph{Step 2: the numerator.} Since $t=p+s$ and
$|s|\le\frac12\min(p,q)$, Lemma~\ref{lem:D-taylor} gives
\[
D(t\,\|\,p)=\frac{s^{2}}{2pq}+R,\qquad |R|\le C_{p}|s|^{3}.
\]
Multiply by $n$. From $s=x\sqrt{pq}\,n^{-1/2}$,
\begin{align*}
n\cdot\frac{s^{2}}{2pq}&=n\cdot\frac{x^{2}pq}{2pq\,n}=\frac{x^{2}}{2},\\
n|R|\le C_{p}\,n|s|^{3}
&=C_{p}\,(pq)^{3/2}|x|^{3}\,n\cdot n^{-3/2}
=C_{p}\,(pq)^{3/2}|x|^{3}\,n^{-1/2}.
\end{align*}
So $u=nD(t\,\|\,p)-\frac{x^{2}}{2}=nR$ satisfies
$|u|\le C_{p}(pq)^{3/2}|x|^{3}\,n^{-1/2}$. The coefficient
$C_{p}(pq)^{3/2}|x|^{3}$ is finite, because $x$ is finite, so $u$ is
$n^{-1/2}$-close to $0$. By the consequence of \eqref{eq:transfer}
noted after Lemma~\ref{lem:relative}, $e^{-u}$ is relatively
$n^{-1/2}$-close to $1$. Since
$e^{-nD(t\|p)}=e^{-x^{2}/2}\,e^{-u}$, this says
\[
e^{-nD(t\|p)}\ \text{is relatively } n^{-1/2}\text{-close to}\ e^{-x^{2}/2}.
\]

\emph{Step 3: the denominator.} Expanding,
\[
t(1-t)=(p+s)(q-s)=pq+s(q-p)-s^{2},
\qquad\text{so}\qquad
\frac{t(1-t)}{pq}=1+\frac{s(q-p)}{pq}-\frac{s^{2}}{pq}.
\]
Using $|s|=|x|\sqrt{pq}\,n^{-1/2}$ and $|q-p|\le 1$,
\[
\Bigl|\frac{s(q-p)}{pq}\Bigr|\le\frac{|x|}{\sqrt{pq}}\,n^{-1/2},
\qquad
\frac{s^{2}}{pq}=\frac{x^{2}}{n}\le x^{2}\,n^{-1/2}.
\]
Both coefficients are finite, so $r=t(1-t)/(pq)$ is $n^{-1/2}$-close
to $1$. Taking square roots preserves this: since $r>0$,
\[
\bigl|\sqrt r-1\bigr|=\frac{|r-1|}{\sqrt r+1}\le|r-1|.
\]
Hence $\sqrt{2\pi n\,t(1-t)}=\sqrt{2\pi npq}\,\sqrt r$ is relatively
$n^{-1/2}$-close to $\sqrt{2\pi npq}$.

\emph{Step 4: conclusion.} Since $n^{-1}\le n^{-1/2}$, the closeness in
\eqref{eq:dML-D} also holds with $n^{-1/2}$ in place of $n^{-1}$. By
Steps 2 and 3 and Lemma~\ref{lem:relative} (quotients, then
transitivity), $b(k)$ is relatively $n^{-1/2}$-close to
\[
\frac{e^{-x^{2}/2}}{\sqrt{2\pi npq}}=\frac{\varphi(x)}{\sqrt{npq}}.
\qedhere
\]
\end{proof}

\begin{example}
For $\omega$ fair flips, $\omega/2$ is an omnific integer and $x=0$,
so the probability of exactly $\omega/2$ heads is relatively
$\omega^{-1/2}$-close to $\sqrt{2/(\pi\omega)}$. This is an
infinitesimal of order $\varepsilon^{1/2}$, compared with
$2^{-\omega}$ for no heads.
\end{example}

\begin{corollary}[Local uniformity]\label{cor:local-uniform}
Under the hypotheses of Theorem~\ref{thm:dML}, if $k'$ is omnific and
$k'-k$ is finite, then $b(k')$ is relatively $n^{-1/2}$-close to
$b(k)$. Consequently, for any finite $m$ and any $m$ consecutive omnific
integers $k,k+1,\dots,k+m-1$ near the mean (that is, with $x$ finite),
the conditional probability of each one, given that $X$ lies among
them, is relatively $n^{-1/2}$-close to $1/m$.
\end{corollary}

\begin{proof}
With $x'=(k'-np)/\sqrt{npq}$, the difference $x'-x=(k'-k)/\sqrt{npq}$
is $n^{-1/2}$-close to $0$. Hence $x'^{2}-x^{2}=(x'-x)(x'+x)$ is too, and
$\varphi(x')/\varphi(x)$ is relatively $n^{-1/2}$-close to $1$. Apply
Theorem~\ref{thm:dML} to $k$ and $k'$. The second statement involves
only a finite sum.
\end{proof}

Thus, viewed through a ``microscope'' of finite width, the binomial
distribution with infinitely many trials looks like the uniform
distribution of Section 2. The bell shape appears only at the scale
$\sqrt{n}$.

\begin{remark}\label{rem:normalization}
For finite $k$ the values $b(k)$ need no hypothesis at all. For infinite
$k$ we have used (F) but never the normalization $\sum_{k}b(k)=1$.
Making sense of that sum, or of expectations $\sum_{k}g(k)b(k)$ that
would extend Theorem~\ref{thm:uniform-close-exp}, requires a theory of summation over
uncountable sets of omnific integers that we do not attempt here.
There is also the question of a sample space: $\omega$ independent
trials suggests the space $\{0,1\}^{\mathbb{N}}$, with the number of
successes being the magnum of the set of successful trials. However,
a typical subset of $\mathbb{N}$ is not eventually periodic, and is not
known to have a magnum even in the broader sense of Lynch and Mackey
(Remark~\ref{rem:magnum-facts}), so Definition~\ref{def:events} does
not apply directly. For now, $b$ should be read
as a formally defined distribution, just as the urn with $a\omega$
chips is.
\end{remark}

\section{Infinitely Many Draws}

Return to the urn of Section 3, with $K=a\omega$ red and $b\omega$ white
chips, $N=(a+b)\omega$, and $p=a/(a+b)$. Now draw $n$ chips, where $n$
is an omnific integer with both $n$ and $N-n$ infinite. Let
\[
h(k)=\frac{\binom{a\omega}{k}\binom{b\omega}{n-k}}{\binom{N}{n}},\qquad
\tau=\frac nN,\qquad
\sigma^{2}=npq\,(1-\tau)=npq\,\frac{N-n}{N}.
\]
The quantity $\tau$ is the \emph{sampling fraction}, and $1-\tau$ is
(up to a factor relatively $\varepsilon$-close to $1$) the classical
finite-population correction $(N-n)/(N-1)$. Note that
$\sigma^{2}\ge \frac{pq}{2}\min(n,N-n)$, so $\sigma$ is infinite.

The proof of the next theorem uses the analogue of
Lemma~\ref{lem:D-taylor} for the entropy function $\mathrm{H}$. It
follows directly from that lemma.

\begin{lemma}\label{lem:H-taylor}
Let $0<p<1$ be real, $q=1-p$, and let $C_{p}$ be as in
Lemma~\ref{lem:D-taylor}. Then for every surreal $s$ with
$|s|\le\frac12\min(p,q)$,
\[
\mathrm{H}(p+s)=\mathrm{H}(p)+s\log\frac qp-\frac{s^{2}}{2pq}-R(s),
\qquad |R(s)|\le C_{p}|s|^{3}.
\]
\end{lemma}

\begin{proof}
Put $t=p+s$, so $1-t=q-s$. From the definitions of $\mathrm{H}$ and $D$,
\[
D(t\,\|\,p)=-\mathrm{H}(t)-t\log p-(1-t)\log q.
\]
Here
\[
\begin{aligned}
t\log p+(1-t)\log q&=(p+s)\log p+(q-s)\log q\\
&=\bigl(p\log p+q\log q\bigr)+s\bigl(\log p-\log q\bigr)\\
&=-\mathrm{H}(p)-s\log\frac qp,
\end{aligned}
\]
so
\[
D(p+s\,\|\,p)=-\mathrm{H}(p+s)-\Bigl(-\mathrm{H}(p)-s\log\frac qp\Bigr)
=-\mathrm{H}(p+s)+\mathrm{H}(p)+s\log\frac qp,
\]
that is,
\[
\mathrm{H}(p+s)=\mathrm{H}(p)+s\log\frac qp-D(p+s\,\|\,p).
\]
By Lemma~\ref{lem:D-taylor},
$D(p+s\,\|\,p)=\frac{s^{2}}{2pq}+R(s)$ with $|R(s)|\le C_{p}|s|^{3}$.
\end{proof}

\begin{theorem}\label{thm:hyper-local}
Assume (F). If $k$ is omnific and $x=(k-np)/\sigma$ is finite, then $h(k)$
is relatively $\sigma^{-1}$-close to $\varphi(x)/\sigma$.
\end{theorem}

\begin{proof}
Throughout, $K=a\omega=pN$, so that $b\omega=N-K=qN$, and
$\delta=k-np=x\sigma$. We use repeatedly that
\begin{equation}\label{eq:sigma-bounds}
\sigma^{2}\le npq
\qquad\text{and}\qquad
\sigma^{2}\le (N-n)pq,
\end{equation}
which hold because $\sigma^{2}=npq\cdot\frac{N-n}{N}$ and
$\sigma^{2}=(N-n)pq\cdot\frac nN$, with $\frac{N-n}{N}\le1$ and
$\frac nN\le1$. Equivalently, $\frac1n\le\frac{pq}{\sigma^{2}}$ and
$\frac1{N-n}\le\frac{pq}{\sigma^{2}}$. Since $\sigma$ is infinite, in
particular $\sigma\ge1$, so $\sigma^{-2}\le\sigma^{-1}$.

\emph{Step 1: rearranging the binomial coefficients.} Writing each
binomial coefficient as a quotient of factorials,
\[
h(k)=\frac{\binom{K}{k}\binom{N-K}{n-k}}{\binom Nn}
=\frac{K!\,(N-K)!\,n!\,(N-n)!}{N!\,k!\,(K-k)!\,(n-k)!\,(N-K-n+k)!}
=\frac{\binom nk\binom{N-n}{K-k}}{\binom NK}.
\]
The second form is more convenient, because its top and bottom
indices are close to the fixed proportion $p$, as we now check.

\emph{Step 2: the indices are infinite.} Put
\[
t_{1}=\frac kn=p+s_{1},\quad s_{1}=\frac{\delta}{n},
\qquad
t_{2}=\frac{K-k}{N-n}=p+s_{2},\quad s_{2}=-\frac{\delta}{N-n},
\]
where the second formula uses $K-k=pN-np-\delta=p(N-n)-\delta$. By
\eqref{eq:sigma-bounds},
\[
|s_{1}|=\frac{|x|\sigma}{n}\le|x|\sigma\cdot\frac{pq}{\sigma^{2}}
=\frac{pq\,|x|}{\sigma},
\qquad\text{and likewise}\qquad
|s_{2}|\le\frac{pq\,|x|}{\sigma}.
\]
Since $x$ is finite and $\sigma$ is infinite, $s_{1}$ and $s_{2}$ are
infinitesimal; in particular $|s_{i}|\le\frac12\min(p,q)$. Therefore
\begin{align*}
k&=n(p+s_{1})\ge\tfrac12 np,\\
n-k&=n(q-s_{1})\ge\tfrac12 nq,\\
K-k&=(N-n)(p+s_{2})\ge\tfrac12(N-n)p,\\
(N-n)-(K-k)&=(N-n)(q-s_{2})\ge\tfrac12(N-n)q,
\end{align*}
and $K=pN$, $N-K=qN$. All of these are infinite, so
Lemma~\ref{lem:stirling-binom} applies to each of
$\binom nk$, $\binom{N-n}{K-k}$ and $\binom NK$.

\emph{Step 3: the error terms.} The error terms from
Lemma~\ref{lem:stirling-binom} are $1/\min(k,n-k)$,
$1/\min\bigl(K-k,(N-n)-(K-k)\bigr)$ and $1/\min(K,N-K)$. Since
$\frac{1}{\min(\alpha,\beta)}\le\frac1\alpha+\frac1\beta$ for positive $\alpha,\beta$, they satisfy
\begin{align*}
\frac1k+\frac1{n-k}&\le\frac{2}{np}+\frac{2}{nq}=\frac{2}{npq}\le\frac{2}{\sigma^{2}},\\
\frac1{K-k}+\frac1{(N-n)-(K-k)}&\le\frac{2}{(N-n)p}+\frac{2}{(N-n)q}=\frac{2}{(N-n)pq}\le\frac{2}{\sigma^{2}},\\
\frac1K+\frac1{N-K}&=\frac{1}{Np}+\frac{1}{Nq}=\frac{1}{Npq}\le\frac{1}{npq}\le\frac{1}{\sigma^{2}}.
\end{align*}
In each line, the first inequality uses the lower bounds from Step~2
(and $K=pN$, $N-K=qN$ in the last line), and the equalities use
$\frac1p+\frac1q=\frac{p+q}{pq}=\frac1{pq}$. The final inequalities
are \eqref{eq:sigma-bounds} with reciprocals taken: $\sigma^{2}\le npq$
gives $\frac{1}{npq}\le\frac{1}{\sigma^{2}}$, and
$\sigma^{2}\le(N-n)pq$ gives $\frac{1}{(N-n)pq}\le\frac{1}{\sigma^{2}}$.
In the last line we also used $n\le N$. So all three binomial coefficients are relatively $\sigma^{-2}$-close
to their Stirling approximations, and by Lemma~\ref{lem:relative}
(products and quotients) $h(k)$ is relatively $\sigma^{-2}$-close to
\[
\frac{e^{n\mathrm{H}(t_{1})}}{\sqrt{2\pi n\,t_{1}(1-t_{1})}}
\cdot
\frac{e^{(N-n)\mathrm{H}(t_{2})}}{\sqrt{2\pi(N-n)\,t_{2}(1-t_{2})}}
\cdot
\frac{\sqrt{2\pi Npq}}{e^{N\mathrm{H}(p)}}
=e^{E}\cdot P,
\]
where
\begin{align*}
E&=n\mathrm{H}(t_{1})+(N-n)\mathrm{H}(t_{2})-N\mathrm{H}(p),\\
P&=\frac{1}{\sqrt{2\pi}}
\Bigl(\frac{n(N-n)}{N}\cdot\frac{t_{1}(1-t_{1})\,t_{2}(1-t_{2})}{pq}\Bigr)^{-1/2}.
\end{align*}
We treat $E$ and $P$ separately.

\emph{Step 4: the exponent.} By Lemma~\ref{lem:H-taylor} applied with
$s=s_{1}$ and $s=s_{2}$,
\begin{align*}
n\mathrm{H}(t_{1})&=n\mathrm{H}(p)+ns_{1}\log\frac qp-\frac{ns_{1}^{2}}{2pq}-nR(s_{1}),\\
(N-n)\mathrm{H}(t_{2})&=(N-n)\mathrm{H}(p)+(N-n)s_{2}\log\frac qp-\frac{(N-n)s_{2}^{2}}{2pq}-(N-n)R(s_{2}).
\end{align*}
Add these and subtract $N\mathrm{H}(p)$, collecting terms by degree in $s_i$.
\begin{itemize}
\item The constant terms cancel: $n\mathrm{H}(p)+(N-n)\mathrm{H}(p)-N\mathrm{H}(p)=0$.
\item The linear terms cancel: $ns_{1}+(N-n)s_{2}=\delta-\delta=0$.
\item The quadratic terms give exactly $-x^{2}/2$:
\begin{align*}
-\frac{1}{2pq}\Bigl(ns_{1}^{2}+(N-n)s_{2}^{2}\Bigr)
&=-\frac{\delta^{2}}{2pq}\Bigl(\frac1n+\frac1{N-n}\Bigr)
=-\frac{\delta^{2}}{2pq}\cdot\frac{N}{n(N-n)}\\
&=-\frac{\delta^{2}}{2\sigma^{2}}
=-\frac{x^{2}}{2},
\end{align*}
using $\sigma^{2}=pq\,n(N-n)/N$ and $\delta=x\sigma$.
\item The remaining term is $u=-nR(s_{1})-(N-n)R(s_{2})$. By
Lemma~\ref{lem:H-taylor}, and since $s_{1}=\delta/n$ and
$s_{2}=-\delta/(N-n)$,
\[
n|R(s_{1})|\le C_{p}\,n|s_{1}|^{3}=C_{p}\frac{|\delta|^{3}}{n^{2}},
\quad
(N-n)|R(s_{2})|\le C_{p}\,(N-n)|s_{2}|^{3}=C_{p}\frac{|\delta|^{3}}{(N-n)^{2}}.
\]
Squaring the reciprocal forms of \eqref{eq:sigma-bounds},
$\frac1n\le\frac{pq}{\sigma^{2}}$ and
$\frac1{N-n}\le\frac{pq}{\sigma^{2}}$, gives
\[
\frac1{n^{2}}\le\frac{(pq)^{2}}{\sigma^{4}}
\qquad\text{and}\qquad
\frac1{(N-n)^{2}}\le\frac{(pq)^{2}}{\sigma^{4}}.
\]
Together with $|\delta|^{3}=|x|^{3}\sigma^{3}$, this gives
\[
|u|\le C_{p}\Bigl(\frac{|\delta|^{3}}{n^{2}}+\frac{|\delta|^{3}}{(N-n)^{2}}\Bigr)
\le 2C_{p}\,|x|^{3}\sigma^{3}\cdot\frac{(pq)^{2}}{\sigma^{4}}
=2C_{p}(pq)^{2}|x|^{3}\,\sigma^{-1}.
\]
\end{itemize}
So $E=-x^{2}/2+u$, where the coefficient $2C_{p}(pq)^{2}|x|^{3}$ is
finite, i.e. $u$ is $\sigma^{-1}$-close to $0$. By the consequence of
\eqref{eq:transfer} noted after Lemma~\ref{lem:relative}, $e^{u}$ is
relatively $\sigma^{-1}$-close to $1$, and hence
\[
e^{E}=e^{-x^{2}/2}e^{u}\ \text{is relatively } \sigma^{-1}\text{-close to}\ e^{-x^{2}/2}.
\]

\emph{Step 5: the prefactor.} Since
$\frac{n(N-n)}{N}=\frac{\sigma^{2}}{pq}$,
\[
P=\frac{1}{\sqrt{2\pi}\,\sigma}\,\bigl(r_{1}r_{2}\bigr)^{-1/2},
\qquad
r_{i}=\frac{t_{i}(1-t_{i})}{pq}=1+\frac{s_{i}(q-p)}{pq}-\frac{s_{i}^{2}}{pq},
\]
where the last formula comes from $t_{i}(1-t_{i})=(p+s_{i})(q-s_{i})$, as
in the proof of Theorem~\ref{thm:dML}. By the bound on $|s_{i}|$ in
Step 2 and $|q-p|\le1$,
\[
\Bigl|\frac{s_{i}(q-p)}{pq}\Bigr|\le\frac{|x|}{\sigma},
\qquad
\frac{s_{i}^{2}}{pq}\le\frac{pq\,x^{2}}{\sigma^{2}}\le\frac{pq\,x^{2}}{\sigma},
\]
so each $r_{i}$ is $\sigma^{-1}$-close to $1$. By
Lemma~\ref{lem:relative}, so is $r_{1}r_{2}$; by the square-root
inequality in the proof of Theorem~\ref{thm:dML}, so is
$\sqrt{r_{1}r_{2}}$; and by Lemma~\ref{lem:relative} again, so is
$(r_{1}r_{2})^{-1/2}$. Hence $P$ is relatively $\sigma^{-1}$-close to
$1/(\sqrt{2\pi}\,\sigma)$.

\emph{Step 6: conclusion.} By Step 3, $h(k)$ is relatively
$\sigma^{-2}$-close, and hence $\sigma^{-1}$-close, to $e^{E}P$. By
Steps 4 and 5 and Lemma~\ref{lem:relative}, $e^{E}P$ is relatively
$\sigma^{-1}$-close to
\[
\frac{e^{-x^{2}/2}}{\sqrt{2\pi}\,\sigma}=\frac{\varphi(x)}{\sigma},
\]
and transitivity finishes the proof.
\end{proof}

The binomial and hypergeometric distributions are therefore both
locally normal, with variances $npq$ and $npq(1-\tau)$. Whether they
agree depends only on $\tau$.

\begin{theorem}\label{thm:dichotomy}
Assume (F), and let $n$ be an omnific integer with $n$ and $N-n$ infinite.
\begin{enumerate}
\item If $\tau$ is infinitesimal, then for every omnific $k$ with
$(k-np)/\sqrt{npq}$ finite, $h(k)$ is relatively
$\max\{\tau,n^{-1/2}\}$-close to $b(k;n,p)$.
\item If $\tau$ is not infinitesimal, let $k$ be the omnific integer
with $k\le np<k+1$ \cite{Conway}. Then $h(k)/b(k;n,p)$ is relatively
$\sigma^{-1}$-close to $(1-\tau)^{-1/2}$. Since
$(1-\tau)^{-1/2}-1\ge\tau/2$, the ratio $h(k)/b(k;n,p)$ is not
infinitesimally close to $1$.
\end{enumerate}
\end{theorem}

\begin{proof}
\emph{Step 1: a formula for the ratio.} Let
\[
x=\frac{k-np}{\sqrt{npq}}
\qquad\text{and}\qquad
x_{h}=\frac{k-np}{\sigma}.
\]
Since $\sigma^{2}=npq(1-\tau)$, we have $\sigma=\sqrt{npq}\,\sqrt{1-\tau}$
and hence $x_{h}=x/\sqrt{1-\tau}$. Suppose that $x$ and $x_{h}$ are both
finite (we check this in each case below). Then
Theorem~\ref{thm:dML} gives that $b(k)$ is relatively
$n^{-1/2}$-close to $\varphi(x)/\sqrt{npq}$, and
Theorem~\ref{thm:hyper-local} gives that $h(k)$ is relatively
$\sigma^{-1}$-close to $\varphi(x_{h})/\sigma$. To compare these, note
that $\sigma^{2}\le npq$ by \eqref{eq:sigma-bounds}, so
\[
n^{-1/2}=\frac{\sqrt{pq}}{\sqrt{npq}}\le\sqrt{pq}\,\sigma^{-1},
\]
and the first closeness also holds with $\sigma^{-1}$ in place of
$n^{-1/2}$. By Lemma~\ref{lem:relative} (quotients), $h(k)/b(k)$ is
relatively $\sigma^{-1}$-close to
\[
Q=\frac{\varphi(x_{h})/\sigma}{\varphi(x)/\sqrt{npq}}
=\frac{\sqrt{npq}}{\sigma}\cdot\frac{\varphi(x_{h})}{\varphi(x)}.
\]
Here $\sqrt{npq}/\sigma=(1-\tau)^{-1/2}$, and
\[
\frac{\varphi(x_{h})}{\varphi(x)}
=\exp\Bigl(-\frac{x^{2}}{2(1-\tau)}+\frac{x^{2}}{2}\Bigr)
=\exp\Bigl(-\frac{x^{2}}{2}\cdot\frac{\tau}{1-\tau}\Bigr)
=\exp\Bigl(-\frac{\tau x_{h}^{2}}{2}\Bigr).
\]
So
\begin{equation}\label{eq:ratio-Q}
\frac{h(k)}{b(k)}\ \text{is relatively } \sigma^{-1}\text{-close to}\
Q=(1-\tau)^{-1/2}\,e^{-v},
\qquad v=\frac{\tau x_{h}^{2}}{2}=\frac{x^{2}}{2}\cdot\frac{\tau}{1-\tau}.
\end{equation}
We will also use the identity, valid for $0\le\tau<1$,
\begin{equation}\label{eq:tau-identity}
(1-\tau)^{-1/2}-1=\frac{1-\sqrt{1-\tau}}{\sqrt{1-\tau}}
=\frac{\tau}{\sqrt{1-\tau}\,\bigl(1+\sqrt{1-\tau}\bigr)},
\end{equation}
obtained by multiplying the numerator and denominator by
$1+\sqrt{1-\tau}$.

\emph{Step 2: case (1).} Suppose $\tau$ is infinitesimal and $x$ is
finite. Then $\tau\le\frac12$, so $1-\tau\ge\frac12$, and
\begin{itemize}
\item $x_{h}=x/\sqrt{1-\tau}$ satisfies $|x_{h}|\le\sqrt2\,|x|$, so it
is finite and Step 1 applies;
\item $\sigma^{2}=npq(1-\tau)\ge\frac12 npq$, so
$\sigma^{-1}\le\sqrt{2/(pq)}\;n^{-1/2}$. Thus \eqref{eq:ratio-Q} holds
with $n^{-1/2}$ in place of $\sigma^{-1}$.
\end{itemize}
Next we show that $Q$ is relatively $\tau$-close to $1$, one factor at
a time.
\begin{itemize}
\item By \eqref{eq:tau-identity}, using $\sqrt{1-\tau}\ge1/\sqrt2$ and
$1+\sqrt{1-\tau}\ge1$,
\[
0\le(1-\tau)^{-1/2}-1\le\sqrt2\,\tau,
\]
so $(1-\tau)^{-1/2}$ is $\tau$-close to $1$.
\item Since $1-\tau\ge\frac12$, we have $0\le v\le x^{2}\tau$. The
coefficient $x^{2}$ is finite, so $v$ is $\tau$-close to $0$, and by
the consequence of \eqref{eq:transfer} noted after
Lemma~\ref{lem:relative}, $e^{-v}$ is relatively $\tau$-close to $1$.
\end{itemize}
By Lemma~\ref{lem:relative} (products), $Q$ is relatively
$\tau$-close to $1$.

Finally let $\eta=\max\{\tau,n^{-1/2}\}$. Both closenesses above hold
with $\eta$ in place of $n^{-1/2}$ or $\tau$, so by transitivity
$h(k)/b(k)$ is relatively $\eta$-close to $1$. That is, $h(k)$ is
relatively $\eta$-close to $b(k)$.

\emph{Step 3: case (2).} Suppose $\tau$ is not infinitesimal, and let
$k$ be the omnific integer with $k\le np<k+1$; such a $k$ exists
because every surreal number lies between consecutive omnific
integers \cite{Conway}. Then $|k-np|<1$, so
\[
|x|<\frac{1}{\sqrt{npq}}
\qquad\text{and}\qquad
|x_{h}|<\frac1\sigma.
\]
Both are infinitesimal, since $n$ and $\sigma$ are infinite. In
particular they are finite, so Step 1 applies. Since $0\le\tau\le1$,
\[
0\le v=\frac{\tau x_{h}^{2}}{2}\le\frac{x_{h}^{2}}{2}<\frac{1}{2\sigma^{2}}\le\frac{1}{2\sigma},
\]
using $\sigma\ge1$. So $v$ is $\sigma^{-1}$-close to $0$, and $e^{-v}$
is relatively $\sigma^{-1}$-close to $1$. Hence
$Q=(1-\tau)^{-1/2}e^{-v}$ is relatively $\sigma^{-1}$-close to
$(1-\tau)^{-1/2}$, and by \eqref{eq:ratio-Q} and transitivity,
$h(k)/b(k)$ is relatively $\sigma^{-1}$-close to $(1-\tau)^{-1/2}$.

It remains to show that $h(k)/b(k)$ is not infinitesimally close to
$1$. Put $A=(1-\tau)^{-1/2}$, so $A>0$, and define $\rho$ by
\[
\frac{h(k)}{b(k)}=A(1+\rho),\qquad\text{that is,}\qquad
\rho=\frac{h(k)/b(k)}{A}-1.
\]
Since $h(k)/b(k)$ is relatively $\sigma^{-1}$-close to $A$, the ratio
$\bigl(h(k)/b(k)\bigr)/A$ is $\sigma^{-1}$-close to $1$; so $\rho$ is
$\sigma^{-1}$-close to $0$, and in particular infinitesimal. By
\eqref{eq:tau-identity}, using $\sqrt{1-\tau}\le1$ and
$1+\sqrt{1-\tau}\le2$,
\[
A-1\ge\frac\tau2.
\]
If $A$ is infinite, then $h(k)/b(k)=A(1+\rho)\ge\frac12A$ (as
$|\rho|\le\frac12$) is infinite too, and we are done. Otherwise $A$ is
finite. Expanding the product and subtracting $1$,
\[
\frac{h(k)}{b(k)}-1=A+A\rho-1=(A-1)+A\rho.
\]
Since $A>0$, we have $A\rho\ge-A|\rho|$, and so, using $A-1\ge\tau/2$,
\[
\frac{h(k)}{b(k)}-1\ge(A-1)-A|\rho|\ge\frac\tau2-A|\rho|.
\]
Here $\tau/2$ is not infinitesimal, while $A|\rho|$ is a finite number
times an infinitesimal, hence infinitesimal; in particular
$A|\rho|<\tau/4$. Therefore
\[
\frac{h(k)}{b(k)}-1>\frac\tau4,
\]
a positive number that is not infinitesimal, and $h(k)/b(k)$ is not
infinitesimally close to $1$.
\end{proof}

\begin{example}
Take $a=b=1$, so the urn holds $\omega$ red and $\omega$ white chips.
\begin{itemize}
\item Draw $n=\omega^{1/2}$ chips. Then $\tau=\tfrac12\omega^{-1/2}$ is
infinitesimal, and near its center the distribution of red chips is
indistinguishable from the binomial distribution with $\omega^{1/2}$
trials.
\item Draw $n=\omega$ chips, half the urn. Here $N=2\omega$,
$p=q=\tfrac12$ and $\tau=\tfrac12$, and $N-n=\omega$ is infinite. The
mean is $np=\omega/2$, which is an omnific integer, so for $k=\omega/2$
we have $x=0$ in both Theorem~\ref{thm:dML} and
Theorem~\ref{thm:hyper-local}, and $\varphi(0)=1/\sqrt{2\pi}$.
\begin{itemize}
\item For the hypergeometric distribution,
$\sigma^{2}=npq(1-\tau)=\omega\cdot\tfrac14\cdot\tfrac12=\omega/8$, so
by Theorem~\ref{thm:hyper-local} the probability of exactly
$\omega/2$ red chips is relatively $\sigma^{-1}$-close to
\[
\frac{\varphi(0)}{\sigma}=\frac{1}{\sqrt{2\pi}}\sqrt{\frac{8}{\omega}}
=\frac{2}{\sqrt{\pi\omega}}.
\]
Since $\sigma^{-1}=\sqrt8\,\omega^{-1/2}$, this closeness is the same
as relative $\omega^{-1/2}$-closeness.
\item For the binomial distribution with $\omega$ trials,
$npq=\omega/4$, so by Theorem~\ref{thm:dML} the probability of exactly
$\omega/2$ successes is relatively $\omega^{-1/2}$-close to
\[
\frac{\varphi(0)}{\sqrt{npq}}=\frac{1}{\sqrt{2\pi}}\cdot\frac{2}{\sqrt\omega}
=\sqrt{\frac{2}{\pi\omega}}.
\]
\end{itemize}
The ratio is $\bigl(2/\sqrt{\pi\omega}\bigr)\big/\sqrt{2/(\pi\omega)}=\sqrt2$,
in agreement with Theorem~\ref{thm:dichotomy}(2), where
$(1-\tau)^{-1/2}=\sqrt2$. The hypergeometric probability is larger
because its variance $npq(1-\tau)$ is half the binomial variance $npq$.
\end{itemize}
\end{example}

\begin{remark}
Theorem~\ref{thm:dichotomy} says that the boundary between
``hypergeometric $\approx$ binomial'' and ``hypergeometric
$\not\approx$ binomial'' is not finite versus infinite $n$, but
infinitesimal versus non-infinitesimal sampling fraction $n/N$. The
finite-$n$ result of Section 3 is the case where $n$ is finite (and
hence $\tau$ is infinitesimal). The textbook rule of thumb that a
sample of at most about $5\%$ of the population may be treated as
binomial becomes, in the surreal setting, an exact criterion. (Here
$p$ is a positive real; Section~\ref{sec:rare} shows that the criterion
fails when $p$ is infinitesimal.)
\end{remark}

\section{Rare Successes}\label{sec:rare}

So far the success probability $p$, or the proportion of red chips, has
been a positive real number. In this section it is infinitesimal. This
is the surreal counterpart of the classical ``law of rare events'', and
it has a pleasant feature: in every case treated here, the binomial
coefficients that appear have finite lower index, so they are given by
the polynomial formula \eqref{eq:binom-finite} and Hypothesis (F) is not
needed.

We use one elementary estimate repeatedly. For a surreal $x$ and finite
$k$, write $x^{(k)}=x(x-1)\cdots(x-k+1)$ for the falling factorial, so
that $\binom xk=x^{(k)}/k!$.

\begin{lemma}\label{lem:falling}
If $x$ is infinite and $k$ is finite, then $x^{(k)}$ is relatively
$x^{-1}$-close to $x^{k}$.
\end{lemma}

\begin{proof}
We have
\[
\frac{x^{(k)}}{x^{k}}=\prod_{i=0}^{k-1}\Bigl(1-\frac ix\Bigr).
\]
Each factor is $x^{-1}$-close to $1$, since $|(1-i/x)-1|=i\cdot x^{-1}$
and $i$ is finite. As $x^{-1}$ is infinitesimal, applying
Lemma~\ref{lem:relative} (products) $k-1$ times shows that the product
of these finitely many factors is relatively $x^{-1}$-close to $1$.
\end{proof}

\subsection{The Poisson limit}

\begin{theorem}[Surreal Poisson limit]\label{thm:poisson}
Let $n$ be an infinite omnific integer and let $0<p<1$ with $\lambda=np$
finite (so $p$ is infinitesimal). Then for every finite $k\ge0$,
$b(k;n,p)$ is relatively $n^{-1}$-close to
\[
e^{-\lambda}\frac{\lambda^{k}}{k!}.
\]
No hypothesis on factorials is needed.
\end{theorem}

\begin{proof}
Since $k$ is finite, $\binom nk=n^{(k)}/k!$ by \eqref{eq:binom-finite}.
Using $n^{k}p^{k}=\lambda^{k}$, we can write
\[
b(k)=\binom nk p^{k}q^{n-k}
=\frac{n^{(k)}}{n^{k}}\cdot\frac{\lambda^{k}}{k!}\cdot q^{-k}\cdot q^{n},
\]
and we treat the three factors other than $\lambda^{k}/k!$ in turn.
Note first that $p=\lambda\cdot n^{-1}$ with $\lambda$ finite, so
anything that is $p$-close is also $n^{-1}$-close.
\begin{itemize}
\item By Lemma~\ref{lem:falling}, $n^{(k)}/n^{k}$ is relatively
$n^{-1}$-close to $1$.
\item The number $q=1-p$ is $p$-close to $1$. By
Lemma~\ref{lem:relative} (quotients, then products, $k$ times),
$q^{-k}$ is relatively $p$-close, hence $n^{-1}$-close, to $1$.
\item Define $u=n\log q+\lambda$, so that $n\log q=-\lambda+u$ holds
by definition; the content is that $u$ is small. Since $\lambda=np$,
\[
u=n\log(1-p)+np=n\bigl(\log(1-p)+p\bigr).
\]
Since $p$ is infinitesimal, $|p|\le\frac12$, and the second
inequality in \eqref{eq:transfer} with $t=-p$ gives
$|\log(1-p)+p|\le p^{2}$. Multiplying by $n$, and using $p=\lambda/n$,
\[
|u|\le np^{2}=n\cdot\frac{\lambda^{2}}{n^{2}}=\lambda^{2}\cdot n^{-1}.
\]
The coefficient $\lambda^{2}$ is finite, so $u$ is $n^{-1}$-close to
$0$. Hence $q^{n}=e^{n\log q}=e^{-\lambda}e^{u}$, and $e^{u}$ is
relatively $n^{-1}$-close to $1$ by the consequence of
\eqref{eq:transfer} noted after Lemma~\ref{lem:relative}. So $q^{n}$ is
relatively $n^{-1}$-close to $e^{-\lambda}$.
\end{itemize}
Combining the three items with Lemma~\ref{lem:relative} (products)
gives the theorem.
\end{proof}

\begin{example}\label{ex:poisson-omega}
Perform $\omega$ independent trials, each with success probability
$\varepsilon$. Then $\lambda=1$, and for every finite $k$ the
probability of exactly $k$ successes is relatively $\varepsilon$-close
to $e^{-1}/k!$. In particular, the probability of no successes,
$(1-\varepsilon)^{\omega}$, is relatively $\varepsilon$-close to
$e^{-1}$.
\end{example}

\begin{remark}\label{rem:geometric}
The same estimates shed light on the nap of the Introduction. Suppose
that in each minute the man wakes, independently, with probability
$\varepsilon$, so that the minute in which he wakes has the geometric
probability mass function
\[
g(k)=(1-\varepsilon)^{k-1}\varepsilon,\qquad k=1,2,3,\dots.
\]
For finite $k$, $(1-\varepsilon)^{k-1}$ is relatively
$\varepsilon$-close to $1$ (as in the second item of the proof above),
so $g(k)$ is relatively $\varepsilon$-close to $\varepsilon$. On every
finite stretch of time, then, this memoryless model is
indistinguishable from the uniform model $f(k)=\varepsilon$ of the
Introduction. The two separate only on the scale of $\omega$: if $c$ is
a positive real, then $c\omega$ is an omnific integer, and the argument
of the third item gives that $(1-\varepsilon)^{c\omega-1}$ is
relatively $\varepsilon$-close to $e^{-c}$, so $g(c\omega)$ is
relatively $\varepsilon$-close to $e^{-c}\varepsilon$.

The two models differ in how they arise. The geometric model comes
with a mechanism: it is produced by a sequence of independent trials,
one per minute, each with the same waking probability $\varepsilon$.
In particular it is memoryless: given that the man is still asleep
after $k-1$ minutes, the chance that he wakes in minute $k$ is
$\varepsilon$, whatever $k$ is. The uniform model, by contrast, simply
assigns the probabilities. Its conditional waking probability is
\[
\frac{f(k)}{P(X\ge k)}=\frac{\varepsilon}{1-(k-1)\varepsilon},
\]
which is infinitely close to $\varepsilon$ for finite $k$ but grows as
$k$ approaches the scale of $\omega$. (Here $P(X\ge k)$ is computed as
$1-\sum_{j<k}f(j)$, or equivalently as the magnum of
$\{k,k+1,\dots\}$, namely $\omega-k+1$, divided by $\omega$.) So the
uniform model is exactly uniform but not memoryless, while the
geometric model is memoryless but uniform only locally, just like the
binomial distribution in Corollary~\ref{cor:local-uniform}.
\end{remark}

\subsection{Urns with few red chips}

Now consider an urn with $r$ red chips and $W$ white chips,
$N=r+W$ in all, with $N$ infinite. We draw $n$ chips without
replacement and let $h(k)$ be the probability that exactly $k$ of them
are red. We first record the classical symmetry of the hypergeometric
distribution in a form that makes sense when $N$ is surreal.

\begin{lemma}\label{lem:hyper-symmetry}
Let $r$, $n$ and $k$ be finite with $0\le k\le\min(r,n)$, and let $N$
be infinite. Then
\[
\frac{\binom rk\binom{N-r}{n-k}}{\binom Nn}
=\frac{\binom nk\binom{N-n}{r-k}}{\binom Nr}.
\]
\end{lemma}

\begin{proof}
By \eqref{eq:binom-finite}, each side is a rational function of $N$
with real coefficients, and the denominators $\binom Nn$ and $\binom Nr$
are nonzero for infinite $N$. As written, however, the two sides are
\emph{different} rational functions: the left side is built from the
falling factorials $(N-r)^{(n-k)}$ and $N^{(n)}$, the right side from
$(N-n)^{(r-k)}$ and $N^{(r)}$.

For every integer $N\ge r+n$, the two sides are equal: writing all the
binomial coefficients with ordinary factorials, both become
\[
\frac{r!\,(N-r)!\,n!\,(N-n)!}{N!\,k!\,(r-k)!\,(n-k)!\,(N-r-n+k)!}.
\]
This factorial expression only makes sense for integer $N$, so it
cannot be used directly when $N$ is infinite. Instead we argue as
follows. Write the left side as $P_{1}(N)/Q_{1}(N)$ and the right side
as $P_{2}(N)/Q_{2}(N)$, with $P_{i},Q_{i}$ polynomials with real
coefficients. Clearing denominators, the polynomial
\[
\Delta(N)=P_{1}(N)Q_{2}(N)-P_{2}(N)Q_{1}(N)
\]
vanishes at every integer $N\ge r+n$. A nonzero polynomial has only
finitely many roots, so $\Delta$ is the zero polynomial, and hence $\Delta(N)=0$
for every surreal $N$ as well. Dividing by $Q_{1}(N)Q_{2}(N)$, which
is nonzero for infinite $N$, gives the identity.
\end{proof}

When $r$ is finite, the right-hand side of Lemma~\ref{lem:hyper-symmetry}
involves only binomial coefficients with finite lower index, even if $n$
is infinite. The left-hand side does not: when $n$ is infinite,
$\binom{N-r}{n-k}$ has an infinite lower index and requires (F).
Accordingly, when $r$ is finite we \emph{define} $h(k)$ by the
right-hand side,
\begin{equation}\label{eq:h-few-red}
h(k)=\frac{\binom nk\binom{N-n}{r-k}}{\binom Nr},\qquad 0\le k\le r,
\end{equation}
which agrees with the left-hand side whenever (F) is assumed (by the
factorial computation in Step~1 of the proof of
Theorem~\ref{thm:hyper-local}) and makes sense without it.

\begin{theorem}[Finitely many red chips, finitely many draws]\label{thm:few-red}
Let $r$ and $n$ be finite and $N$ infinite. For
$0\le k\le\min(r,n)$, $h(k)$ is relatively $N^{-1}$-close to
\[
\frac{r^{(k)}\,n^{(k)}}{k!\,N^{k}},
\]
and $h(k)=0$ for $k>\min(r,n)$.
\end{theorem}

\begin{proof}
For $k>r$ there are not enough red chips, and for $k>n$ not enough
draws, so $h(k)=0$. Now let $0\le k\le\min(r,n)$ and use
\eqref{eq:h-few-red}. All three binomial coefficients have finite lower
index:
\begin{itemize}
\item $\binom nk=n^{(k)}/k!$ is a finite real number.
\item $\binom{N-n}{r-k}=(N-n)^{(r-k)}/(r-k)!$. Since $N-n$ is infinite,
Lemma~\ref{lem:falling} shows that $(N-n)^{(r-k)}$ is relatively
$(N-n)^{-1}$-close to $(N-n)^{r-k}$; and $(N-n)^{-1}\le 2N^{-1}$. Also
$(N-n)/N=1-n/N$ is $N^{-1}$-close to $1$, so by
Lemma~\ref{lem:relative} $(N-n)^{r-k}$ is relatively $N^{-1}$-close to
$N^{r-k}$. Thus $\binom{N-n}{r-k}$ is relatively $N^{-1}$-close to
$N^{r-k}/(r-k)!$.
\item By Lemma~\ref{lem:falling}, $\binom Nr=N^{(r)}/r!$ is relatively
$N^{-1}$-close to $N^{r}/r!$.
\end{itemize}
By Lemma~\ref{lem:relative}, $h(k)$ is relatively $N^{-1}$-close to
\[
\frac{n^{(k)}}{k!}\cdot\frac{N^{r-k}}{(r-k)!}\cdot\frac{r!}{N^{r}}
=\frac{n^{(k)}}{k!}\cdot\frac{r!}{(r-k)!}\cdot\frac{1}{N^{k}}
=\frac{r^{(k)}\,n^{(k)}}{k!\,N^{k}}.
\qedhere
\]
\end{proof}

So the distribution of the number of red chips is graded by powers of
$N^{-1}$: $h(0)$ is infinitely close to $1$, $h(1)$ is of order
$N^{-1}$, $h(2)$ of order $N^{-2}$, and so on. Classically, letting
$N\to\infty$ with $r$ and $n$ fixed collapses this to the statement
that no red chip is drawn; in the surreal setting every level remains
visible. Note also that the leading term is symmetric in $r$ and $n$.

It is natural to compare $h$ with the binomial distribution with $n$
trials and success probability $p=r/N$, the proportion of red chips.

\begin{theorem}\label{thm:few-red-vs-binom}
Let $r$ and $n$ be finite with $r\ge1$, let $N$ be infinite, and let
$p=r/N$.
\begin{enumerate}
\item For $0\le k\le\min(r,n)$, the ratio $h(k)/b(k;n,p)$ is
relatively $N^{-1}$-close to
\[
\frac{r^{(k)}}{r^{k}}=\prod_{i=0}^{k-1}\Bigl(1-\frac ir\Bigr).
\]
In particular, $h(k)$ is relatively $N^{-1}$-close to $b(k;n,p)$ for
$k=0$ and $k=1$.
\item For $2\le k\le\min(r,n)$, the ratio $h(k)/b(k;n,p)$ is
infinitely close to a real number at most $1-\frac1r$, so $h(k)$ and
$b(k;n,p)$ are not relatively infinitesimally close. For
$r<k\le n$, $h(k)=0$ while $b(k;n,p)>0$.
\item Nevertheless, $h$ and $b(\,\cdot\,;n,p)$ are uniformly
$N^{-2}$-close in the sense of Definition~\ref{def:uniform-close}.
\end{enumerate}
\end{theorem}

\begin{proof}
First we find the leading term of $b(k)=b(k;n,p)$ for
$0\le k\le n$. Since $p^{k}=r^{k}N^{-k}$,
\[
b(k)=\binom nk p^{k}(1-p)^{n-k}=\frac{n^{(k)}\,r^{k}}{k!\,N^{k}}\,(1-p)^{n-k}.
\]
Here $1-p$ is $N^{-1}$-close to $1$ (with constant $r$), and $n-k$ is
finite, so by Lemma~\ref{lem:relative} $(1-p)^{n-k}$ is relatively
$N^{-1}$-close to $1$. Hence
\begin{equation}\label{eq:b-few-red}
b(k)\ \text{is relatively } N^{-1}\text{-close to}\
\frac{n^{(k)}\,r^{k}}{k!\,N^{k}}.
\end{equation}

(1) Divide the fraction in the conclusion of Theorem~\ref{thm:few-red} by the 
fraction in \eqref{eq:b-few-red}, using Lemma~\ref{lem:relative} (quotients). For
$k=0$ and $k=1$ the product $\prod_{i<k}(1-i/r)$ is empty or equal to
$1$.

(2) For $2\le k\le\min(r,n)$, the factor with $i=1$ is $1-\frac1r$ and
every other factor lies in $(0,1]$, so the real number $r^{(k)}/r^{k}$
is at most $1-\frac1r<1$. By (1), $h(k)/b(k)$ differs from it by an
infinitesimal. If $r<k\le n$, then $h(k)=0$ by
Theorem~\ref{thm:few-red}, while $b(k)>0$.

(3) There are only finitely many $k$ (namely $0\le k\le n$) at which
$h(k)$ or $b(k)$ is nonzero, so it suffices to show that
$|h(k)-b(k)|$ is $N^{-2}$-close to $0$ for each such $k$.
\begin{itemize}
\item If $k\ge2$, then by Theorem~\ref{thm:few-red} and
\eqref{eq:b-few-red}, $h(k)$ and $b(k)$ are each at most a finite
multiple of $N^{-k}\le N^{-2}$.
\item If $k=1$, both $h(1)$ and $b(1)$ are relatively $N^{-1}$-close to
$L=rn/N$, so $|h(1)-L|$ and $|b(1)-L|$ are at most finite multiples of
$L\cdot N^{-1}=rn\,N^{-2}$, and hence so is $|h(1)-b(1)|$.
\item If $k=0$, we compare directly. By Lemma~\ref{lem:hyper-symmetry},
$h(0)=\binom{N-r}{n}\big/\binom Nn$, so
\[
h(0)=\prod_{i=0}^{n-1}\frac{N-r-i}{N-i}=\prod_{i=0}^{n-1}\Bigl(1-\frac{r}{N-i}\Bigr),
\qquad
b(0)=\prod_{i=0}^{n-1}\Bigl(1-\frac rN\Bigr).
\]
All factors lie in $[0,1]$. For numbers $a_{1},\dots,a_{m}$ and
$b_{1},\dots,b_{m}$ in $[0,1]$, with $m$ finite, we have
\begin{equation}\label{eq:product-diff}
\Bigl|\prod_{i=1}^{m}a_{i}-\prod_{i=1}^{m}b_{i}\Bigr|\le\sum_{i=1}^{m}|a_{i}-b_{i}|.
\end{equation}
To see this, induct on $m$; the case $m=1$ is trivial. Write
$A=\prod_{i<m}a_{i}$ and $B=\prod_{i<m}b_{i}$, which also lie in
$[0,1]$. Adding and subtracting $a_{m}B$,
\[
a_{m}A-b_{m}B=a_{m}(A-B)+(a_{m}-b_{m})B,
\]
so, since $|a_{m}|\le1$ and $|B|\le1$,
\[
|a_{m}A-b_{m}B|\le|A-B|+|a_{m}-b_{m}|
\le\sum_{i<m}|a_{i}-b_{i}|+|a_{m}-b_{m}|
\]
by the induction hypothesis. Now apply \eqref{eq:product-diff} with
$m=n$, $a_{i}=1-\frac{r}{N-i}$ and $b_{i}=1-\frac rN$ (indexing from
$i=0$ to $n-1$). Since $N-i\ge N/2$,
\[
|a_{i}-b_{i}|=\Bigl|\frac{r}{N-i}-\frac rN\Bigr|=\frac{ri}{N(N-i)}\le\frac{2ri}{N^{2}},
\]
and summing over $0\le i\le n-1$ gives
$|h(0)-b(0)|\le\frac{2r}{N^{2}}\sum_{i=0}^{n-1}i=2r\binom n2 N^{-2}$.\qedhere
\end{itemize}
\end{proof}

\begin{remark}
Theorem~\ref{thm:few-red-vs-binom} is a good illustration of why
relative closeness is the right notion when probabilities are
infinitesimal. In the absolute sense of
Definition~\ref{def:uniform-close}, $h$ and the binomial distribution
are as close as one could wish, because every probability involved is
either infinitely close to $0$ or $1$. Relative closeness sees the
difference: an urn with a single red chip can never yield two red
chips, while the binomial model says that it can, with probability of
order $N^{-2}$. It also shows that the criterion of
Theorem~\ref{thm:dichotomy} depends on $p$ being a positive real: here
the sampling fraction $n/N$ is infinitesimal, yet the hypergeometric
and binomial distributions differ.
\end{remark}

Finally, by Lemma~\ref{lem:hyper-symmetry} the roles of red chips and
draws can be exchanged, and this gives the case of finitely many red
chips and infinitely many draws.

\begin{theorem}[Finitely many red chips, infinitely many draws]\label{thm:few-red-many-draws}
Let $r$ be finite, let $n$ and $N-n$ be infinite, and let $\tau=n/N$.
Then for $0\le k\le r$, $h(k)$ is relatively
$\frac{1}{\min(n,\,N-n)}$-close to
\[
\binom rk\tau^{k}(1-\tau)^{r-k}.
\]
That is, the number of red chips drawn is approximately binomial with
$r$ trials and success probability $\tau$.
\end{theorem}

\begin{proof}
Let $\eta=1/\min(n,N-n)$, and note that $N^{-1}\le\eta$. In
\eqref{eq:h-few-red}, each binomial coefficient has finite lower index
and infinite upper index, so by Lemma~\ref{lem:falling}
\begin{itemize}
\item $\binom nk=n^{(k)}/k!$ is relatively $n^{-1}$-close to
$n^{k}/k!$;
\item $\binom{N-n}{r-k}$ is relatively $(N-n)^{-1}$-close to
$(N-n)^{r-k}/(r-k)!$;
\item $\binom Nr$ is relatively $N^{-1}$-close to $N^{r}/r!$.
\end{itemize}
Each of these error terms is at most $\eta$, so by
Lemma~\ref{lem:relative} $h(k)$ is relatively $\eta$-close to
\[
\frac{n^{k}}{k!}\cdot\frac{(N-n)^{r-k}}{(r-k)!}\cdot\frac{r!}{N^{r}}
=\binom rk\Bigl(\frac nN\Bigr)^{k}\Bigl(\frac{N-n}{N}\Bigr)^{r-k}
=\binom rk\tau^{k}(1-\tau)^{r-k}.
\qedhere
\]
\end{proof}

The interpretation is that each of the $r$ red chips lands in the
sample with probability $\tau$, nearly independently of the others.
If $\tau$ is a positive real, this is an honest binomial distribution;
if $\tau$ is infinitesimal, it is again graded, by powers of $\tau$.

\begin{example}\label{ex:five-red}
Number the chips by the natural numbers, and let the chips numbered
$1,2,3,4,5$ be red and the rest white. Then $r=5$, the white chips have
magnum $\omega-5$, and $N=\omega$, so $N^{-1}=\varepsilon$. Note that
$r^{(k)}/k!=\binom5k$.
\begin{itemize}
\item \emph{Finitely many draws.} Draw $n$ chips, $n$ finite. By
Theorem~\ref{thm:few-red}, the probability of exactly $k$ red chips is
relatively $\varepsilon$-close to $\binom5k n^{(k)}\varepsilon^{k}$
for $0\le k\le\min(5,n)$, and is $0$ for $k>5$. Thus
\begin{itemize}
\item no red chip: relatively $\varepsilon$-close to $1$;
\item exactly one: relatively $\varepsilon$-close to $5n\,\varepsilon$;
\item exactly two: relatively $\varepsilon$-close to
$10\,n(n-1)\,\varepsilon^{2}$;
\item all five (if $n\ge5$): relatively $\varepsilon$-close to
$n(n-1)(n-2)(n-3)(n-4)\,\varepsilon^{5}$;
\item six or more: $0$.
\end{itemize}
Compare the binomial distribution with $n$ trials and success
probability $p=5\varepsilon$, the proportion of red chips. By
Theorem~\ref{thm:few-red-vs-binom}, it agrees with the urn for $k=0$
and $k=1$, but for two red chips it gives approximately
$\binom n2 25\,\varepsilon^{2}=\frac{25}{2}n(n-1)\,\varepsilon^{2}$,
which is $\frac54$ times too large; for five red chips it is too large
by the factor $5^{5}/5!=\frac{625}{24}$; and it gives a positive
probability to six or more red chips, which is impossible. Contrast
Example~\ref{ex:multiples-of-5}, where the proportion of red chips is the positive real
$\frac15$ and the binomial approximation is accurate.
\item \emph{Infinitely many draws.} Draw $n=\omega/2$ chips, so that
$N-n=\omega/2$ is infinite and $\tau=\frac12$. By
Theorem~\ref{thm:few-red-many-draws}, the number of red chips drawn is
approximately binomial with $5$ trials and success probability
$\frac12$: the probability of exactly $k$ red chips is relatively
$\varepsilon$-close to $\binom5k2^{-5}$. For instance, the probability
of drawing all five red chips is relatively $\varepsilon$-close to
$\frac1{32}$. Each red chip is drawn with probability $\frac12$,
essentially independently of the others.
\end{itemize}
\end{example}

\subsection{Summary}

The results for an urn with $N$ chips, $N$ infinite, $r$ of them red,
and $n$ draws can be collected as follows.

\begin{center}
\renewcommand{\arraystretch}{1.3}
\begin{tabular}{p{0.24\textwidth}p{0.2\textwidth}p{0.3\textwidth}p{0.12\textwidth}}
\textbf{red chips $r$} & \textbf{draws $n$} & \textbf{behavior} & \textbf{source}\\
\hline
finite & finite & graded by powers of $N^{-1}$ & Thm~\ref{thm:few-red}\\
finite & $n$, $N-n$ infinite & $\mathrm{Bin}(r,n/N)$ & Thm~\ref{thm:few-red-many-draws}\\
$pN$, $p$ positive real & finite & $\mathrm{Bin}(n,p)$, up to $\varepsilon$ & Section 3\\
$pN$, $p$ positive real & $n$, $N-n$ infinite & locally normal; binomial iff $n/N$ infinitesimal & Thms~\ref{thm:hyper-local}, \ref{thm:dichotomy}\\
\hline
\end{tabular}
\end{center}

\noindent The first three rows require no hypothesis on factorials; the
last requires (F). For the binomial distribution itself, with $n$
infinite, Theorem~\ref{thm:dML} covers $p$ a positive real and
Theorem~\ref{thm:poisson} covers $np$ finite.

Two cases are left open. For the binomial distribution with $p$
infinitesimal but $np$ infinite, we expect Theorem~\ref{thm:dML} to
extend with error $(npq)^{-1/2}$, the surreal form of the rule that $np$
and $nq$ should both be large; the proof would require redoing
Lemma~\ref{lem:D-taylor} in the scaled variable $s/p$. For the urn
with $r$ and $n$ both infinite but $rn/N$ finite, we expect a Poisson
limit with mean $rn/N$; here the binomial coefficients have infinite
lower index, so (F) is needed, and the factorial estimates must be
carried out with some care.

\section{Further directions}There is no shortage of things to work on:
\begin{itemize} 
\item A natural factorial on $\mathbf{Oz}$. Is there a canonical solution of Hypothesis (F), for instance from a surreal Gamma function, that removes the choice of coset representatives and takes omnific values throughout? Remark~\ref{rem:F-exists} gives omnific values on the cosets $c\,\omega^{x}+\mathbb{Z}$ ($c,x$ positive reals). For a general infinite omnific integer $k_0$, such as $\omega+\omega^{1/2}$, is there an omnific integer $k_0!$, relatively $(1/k_0)$-close to $S(k_0)$, that is divisible by $k_0(k_0-1)\cdots(k_0-j+1)$ for every finite $j$? Better still, is there a solution for which every binomial coefficient $\binom nk$ is an omnific integer, so that it could count something (cf. magnums)?

\item A sample-space interpretation of $\omega$ trials via magnums of subsets of $\mathbb{N}$ (Remark~\ref{rem:normalization}).

\item Binomial with $p$ infinitesimal and $np$ infinite and the urn with $r, n$ infinite and $rn/N$ finite as described in Section 6.

\item Surreal extensions of other discrete distributions.

\end{itemize}
\section{Acknowledgements}
The author used Claude Opus 5.5 (Anthropic), an AI system, extensively. Claude proposed the study of infinitely many trials and draws, suggested the explicit construction of magnums on eventually periodic sets (Proposition~\ref{prop:magnum}), drafted statements and proofs of several results (including Lemma~\ref{lem:stirling-binom}, Theorems~\ref{thm:dML}, \ref{thm:hyper-local} and~\ref{thm:dichotomy}, the omnific-valued factorials in Remark~\ref{rem:F-exists}, and the results of Section~6), and suggested expository revisions throughout. The author checked all mathematics and takes full responsibility for the content.

No sponsorships or funding supported this research, and the author has no relevant financial or non-financial competing interests to report.
\printbibliography 

\end{document}